\documentclass[11pt]{article}
\usepackage{amsmath,amssymb}
\begin{document}
\newtheorem{thm}{Theorem}[section]
\newtheorem{prop}[thm]{Proposition}
\newenvironment{dfn}{\medskip\refstepcounter{thm}
\noindent{\bf Definition \thesection.\arabic{thm}\ }}{\medskip}
\newenvironment{cond}{\medskip\refstepcounter{thm}
\noindent{\bf Condition \thesection.\arabic{thm}\ }}{\medskip}
\newenvironment{ex}{\medskip\refstepcounter{thm}
\noindent{\bf Example \thesection.\arabic{thm}\ }}{\medskip}
\newenvironment{rlist}{\begin{list}{$({\rm \roman{enumi}})$}
{\usecounter{enumi}
\setlength{\rightmargin}{10pt}
\setlength{\leftmargin}{40pt}
\setlength{\itemsep}{2pt}
\setlength{\parsep}{0pt}
\setlength{\labelwidth}{40pt}}}{\end{list}}
\newenvironment{proof}[1][,]{\medskip\ifcat,#1
\noindent{{\it Proof}:\ }\else\noindent{\it Proof of #1.\ }\fi}
{\hfill$\square$\medskip}
\def\eq#1{{\rm(\ref{#1})}}
\def\Q{{\mathbb Q}}
\def\R{{\mathbb R}}
\def\O{{\mathbb O}}
\def\Z{{\mathbb Z}}
\def\GL{\mathop{\rm GL}}
\def\SL{\mathop{\rm SL}}
\def\SU{\mathop{\rm SU}}
\def\SO{\mathop{\rm SO}}
\def\U{\mathbin{\rm U}}
\def\sech{\mathop{\rm sech}\nolimits}
\def\sn{\mathop{\rm sn}\nolimits}
\def\cn{\mathop{\rm cn}\nolimits}
\def\dn{\mathop{\rm dn}\nolimits}
\def\d{{\rm d}}
\def\w{\wedge}
\def\C{{\mathbb C}}
\def\Re{\mathop{\rm Re}\nolimits}
\def\Im{\mathop{\rm Im}\nolimits}
\title{Constructing Associative 3-folds by Evolution Equations}
\author{Jason Lotay, Christ Church, Oxford}
\date{}
\maketitle

\linespread{1.1}

\section{Introduction}

This paper gives two methods for constructing associative 3-folds in
$\R^7$, based around the fundamental idea of evolution equations,
and uses them to produce examples. It is a generalisation of the
work by Joyce in \cite{Joy2}, \cite{Joy3}, \cite{Joy4} and
\cite{Joy5} on special Lagrangian (SL) 3-folds in $\C^3$. The
methods described involve the use of an affine evolution equation
with affine evolution data and the area of ruled submanifolds.

We begin in $\S$\ref{G2} by introducing the exceptional Lie group
$\text{G}_2$ and its relationship with the geometry of associative
3-folds in $\R^7$.  In $\S$\ref{Joyreview} we review the work by
Joyce in \cite{Joy2}, \cite{Joy3} and \cite{Joy4} on evolution
equation constructions for SL $m$-folds in $\C^m$.  We follow this
in $\S$\ref{firsteq} with a derivation of an evolution equation for
associative 3-folds.

In $\S$\ref{second} we derive an affine evolution equation using
affine evolution data.  This is used on an example of such data to
construct a 14-dimensional family of associative 3-folds. One of the
main results of the paper is an explicit solution of the system of
differential equations generated in a particular case to give a
12-dimensional family of associative 3-folds. Moreover, we find a
straightforward condition which ensures that the associative 3-folds
constructed are closed and diffeomorphic to
$\mathcal{S}^1\times\R^2$, rather than $\R^3$.

In the final section, $\S$\ref{rul}, we define ruled associative
3-folds and derive an evolution equation for them.  This allows us
to characterise a family of ruled associative 3-folds using a pair
of real analytic maps satisfying two partial differential equations.
We finish by giving a means of constructing ruled associative
3-folds $M$ from r-oriented two-sided associative cones $M_0$ such
that $M$ is asymptotically conical to $M_0$ with order $O(r^{-1})$.

\medskip

\noindent\textbf{Acknowledgements} \hspace{4pt} I owe a great deal
of gratitude to my supervisor Dominic Joyce for his help, guidance
and foundational work.  Thanks are also due to the referees for
their work in providing helpful suggestions and corrections.

\section{Introduction to $\text{G}_2$ and Associative 3-folds}
\label{G2}

We give two equivalent definitions of $\text{G}_2$, which relate to
the geometry of $\R^7$ and the octonions respectively.  The first
follows \cite[page 242]{Joy1}.

\begin{dfn}
\label{dfn1} Let $(x_1,\ldots,x_7)$ be coordinates on $\R^7$.  We
shall write $\d{\bf x}_{ij\ldots k}$ for the form $\d x_i\w\d
x_j\w\ldots \w\d x_k$ on $\R^7$.  Define a 3-form $\varphi$ on
$\R^7$ by
\begin{equation}
\begin{split}
\varphi \, = \; &\d{\bf x}_{123}+\d{\bf x}_{145}
+\d{\bf x}_{167}+\d{\bf x}_{246}-\d{\bf x}_{257}
-\d{\bf x}_{347}-\d{\bf x}_{356}.
\label{phi}
\end{split}
\end{equation}
Then $\text{G}_2$ = $\bigl\{\gamma\in\GL(7,\R) \: : \:
 \gamma^*\varphi = \varphi \bigr\}$.
\label{def1}
\end{dfn}

 We note that $\text{G}_2$ is a compact, connected, simply
connected, simple, 14-dimensional Lie group, which preserves the
Euclidean metric and the orientation on $\R^7$. It also preserves
the 4-form $\ast\varphi$ given by
\begin{equation}
\begin{split}
\ast\varphi \, = \; &\d{\bf x}_{4567}+\d{\bf x}_{2367}
+\d{\bf x}_{2345}+\d{\bf x}_{1357}-\d{\bf x}_{1346}
-\d{\bf x}_{1256}-\d{\bf x}_{1247},
\label{starphi}
\end{split}
\end{equation}
where $\varphi$ and $\ast\varphi$ are related by the Hodge star.

The second definition, taken from \cite{HarLaw}, comes from
considering the algebra of the octonions, or Cayley numbers, $\O$.

\begin{dfn}
The group of automorphisms of $\O$ is $\text{G}_2$. \label{def2}
\end{dfn}

Suppose we take the latter definition of $\text{G}_2$ and note that
$x \in$ Im$\,\O$
 if and only if $x^2$ is real but $x$ is not.  Therefore, for all $\gamma\in \text{G}_2$
and for $x\in\O$, $\gamma(x) \in$ Im$\,\O$ $\Leftrightarrow$
$\gamma(x)^2 = \gamma(x^2) \in\R$, $\gamma(x)\notin\R$
$\Leftrightarrow$ $x^2\in\R$, $x\notin\R$ $\Leftrightarrow$ $x \in$
Im$\,\O$. Hence, $\text{G}_2$ is the subgroup of the group of
automorphisms of Im$\,\O\cong\R^7$ preserving the octonionic
multiplication on Im$\,\O$. This multiplication defines a cross
product $\times:\R^7\times\R^7 \rightarrow\R^7$ by
\begin{equation}
\begin{split}
x \times y \, = \, &\frac{1}{2}\,(xy - yx), \label{cross1}
\end{split}
\end{equation}
where the right-hand side is defined by considering $x$ and $y$ as
imaginary octonions.  Note that we can recover the octonionic
multiplication from the cross product and also that the cross
product can be written as follows:
\begin{equation}
\begin{split}
(x \times y)^d \, = \, \varphi_{abc}x^a y^b  g^{cd} \label{cross2}
\end{split}
\end{equation}
using index notation for tensors on $\R^7$, where $g^{cd}$ is the
inverse of the\linebreak Euclidean metric on $\R^7$.  This can be
verified using \eq{phi}, \eq{cross1} and a Cayley multiplication
table for the octonions.  We deduce from \eq{cross2} that
\begin{equation}
\begin{split}
\label{innerphi} \varphi(x,y,z) = g(x\times y,z)
\end{split}
\end{equation}
for $x,y,z\in\R^7$, where $g$ is the Euclidean metric on $\R^7$.

\medskip

For this article, we take manifolds to be smooth and nonsingular
almost everywhere and submanifolds to be immersed, unless otherwise
stated. We define \emph{calibrations} and \emph{calibrated
submanifolds} following the approach in \cite{HarLaw}.

\begin{dfn}
Let $(M,g)$ be a Riemannian manifold.  An \emph{oriented tangent\
$k$-plane} $V$ on $M$ is an oriented $k$-dimensional vector subspace
$V$ of $T_xM$, for some $x$ in $M$. Given an oriented tangent
$k$-plane $V$ on $M$, $g|_V$ is a\linebreak Euclidean metric on $V$
and hence, using $g|_V$ and the orientation on $V$, we have a
natural volume
form, vol$_V$, which is a $k$-form on $V$.  \\
Let $\eta$ be a closed $k$-form on $M$.  Then $\eta$ is a
\emph{calibration} on $M$ if $\eta | _V \leq$ vol$_V$ for all
oriented tangent $k$-planes $V$ on $M$, where $\eta | _V = \alpha
\cdot\text{vol}_V$ for some
$\alpha\in\R$, and so $\eta | _V \leq$ vol$_V$ if $\alpha\leq 1$. \\
Let $N$ be an oriented $k$-dimensional submanifold of $M$.  Then $N$
is a \linebreak\emph{calibrated submanifold} or \emph{$\eta
-$submanifold} if $\eta | _{T_xN} =$ vol$_{T_xN}$ for all $x\in N$.
\label{calibration}
\end{dfn}

 Calibrated submanifolds are \textit{minimal} submanifolds
\cite[Theorem II.4.2]{HarLaw}.  We now define \emph{associative
3-folds}.

\begin{dfn}
Let $N$ be a 3-dimensional submanifold of $\R^7$.  Note that, by
\cite[Theorem IV.1.4]{HarLaw}, $\varphi$ as given by \eq{phi} is a
calibration on $\R^7$.  An oriented 3-plane $V$ in $\R^7$ is
\emph{associative} if $\varphi | _V$ = vol$_V$.  $N$ is an
\emph{associative 3-fold} if $T_xN$ is associative for all $x\in N$,
i.e. if $N$ is a $\varphi -$submanifold.
\end{dfn}

An alternative description of associative 3-planes is given in
\cite{HarLaw} which requires the definition of the \emph{associator}
of three octonions.

\begin{dfn}
\label{assdfn} The \emph{associator} $[x,y,z]$ of $x,y,z\in\O$ is
given by
\begin{equation}
\begin{split}
\label{associator}
[x,y,z] \,  = \, (xy)z - x(yz).
\end{split}
\end{equation}
\end{dfn}
\vspace{-20pt}

\noindent Whereas the commutator measures the extent to which
commutativity fails, the associator gives the degree to which
associativity fails in $\O$.  Note that we can write an alternative
formula, in index notation, for the associator of three vectors
$x,y,z \in \R^7$ using $*\varphi$ and the inverse of the Euclidean
metric $g$ on $\R^7$ as follows:
\begin{equation}
\begin{split}\label{asctr}
\frac{1}{2} [x,y,z]^{e} = (*\varphi)_{abcd}x^ay^bz^cg^{de}.
\end{split}
\end{equation}
\noindent This can be verified using \eq{starphi}, \eq{associator}
and a Cayley multiplication table for $\O$.
 We then have the following result \cite[Corollary IV.1.7]{HarLaw}.

\begin{prop}
\label{ass} Let $V$ be a 3-plane in \emph{Im} $\!\O \cong \!\R^7$
with basis $(x,y,z)$.  Then $V$, with an appropriate orientation, is
associative if and only if $[x,y,z] = 0$.
\end{prop}

\noindent In $\S$\ref{second} we require some properties of the
associator which we state as a proposition taken from
\cite[Proposition IV.B.16]{HarLaw}.

\begin{prop}
\label{assprops} The associator $[x,y,z]$ of $x,y,z\in\O$ is:
\begin{rlist}
\item alternating,
\item imaginary valued,
\item orthogonal to $x,y,z$ and to $[a,b] = ab-ba$ for any subset $\{a,b\}$ of
$\{x,y,z\}$.
\end{rlist}
\end{prop}

\section{Special Lagrangian
$m$-folds in $\C^m$ \cite{Joy2} \cite{Joy3} \cite{Joy4}}
\setcounter{thm}{0}
\label{Joyreview}

We review the work in Joyce's papers \cite{Joy2}, \cite{Joy3} and
\cite{Joy4} on the construction of special Lagrangian (SL) $m$-folds
in $\C^m$ using evolution equations, upon which this paper is based.
 We begin by defining the SL calibration form on $\C^m$ and hence SL $m$-folds.

\begin{dfn}
Let $(z_1,\ldots,z_m)$ be complex coordinates on $\C^m$
with\linebreak complex structure $I$.  Define a metric $g$, a real
2-form $\omega$ and a complex $m$-form $\Omega$ on $\C^m$ by
\begin{align*}
g & =  |dz_1|^2 + \ldots + |dz_m|^2, \\
\omega & =  \frac{i}{2}(dz_1\w d\bar{z}_1 + \ldots + dz_m\w d\bar{z}_m), \\
\Omega & =  dz_1 \w \ldots \w dz_m.
\end{align*}
Let $L$ be a real oriented $m$-dimensional submanifold of $\C^m$.
Then $L$ is a {\it special Lagrangian} (SL) $m$-fold in $\C^m$ with
{\it phase} $e^{i\theta}$ if $L$ is calibrated with respect to the
real $m$-form $\cos\theta\Re\Omega+\sin\theta\Im\Omega$. If the
phase of $L$ is unspecified it is taken to be one so that $L$ is
calibrated with respect to $\Re\Omega$.
\end{dfn}

\vspace{-12pt}

Harvey and Lawson \cite[Corollary III.1.11]{HarLaw} give the
following alternative characterisation of SL $m$-folds.

\begin{prop}
Let $L$ be a real $m$-dimensional submanifold of $\C^m$.  Then $L$
admits an orientation making it into an SL $m$-fold in $\C^m$ with
phase $e^{i\theta}$ if and only if $\omega | _L \equiv 0$ and
$(\sin\theta\Re\Omega - \cos \theta\Im\Omega)|_L \equiv 0$.
\end{prop}

Joyce, in \cite{Joy2}, derives an evolution equation for SL
$m$-folds, the proof of which requires the following result
\cite[Theorem III.5.5]{HarLaw}.

\begin{thm}
Let $P$ be a real analytic $(m-\!1)$-dimensional submanifold of
$\C^m\!$ with $\omega|_P \!\equiv 0$.
 Then there exists a unique SL $m$-fold in $\C^m$ containing $P$.
\end{thm}
The requirement that $P$ be real analytic is due to the fact that
the proof uses the {\it Cartan--K\"ahler Theorem}, which is only
applicable in the real analytic category.  We now give the main
result \cite[Theorem 3.3]{Joy2}.

\begin{thm}
Let $P$ be a compact, orientable, $(m-1)$-dimensional, real analytic
manifold, let $\chi$ be a real analytic nowhere vanishing section of
$\Lambda^{m-1}TP$ and let $\psi:P \rightarrow\C^m$ be a real
analytic embedding (immersion) such that $\psi^*(\omega)\equiv 0$ on
$P$. Then there exist $\epsilon > 0$ and a unique family
$\{\psi_t:t\in(-\epsilon,\epsilon)\}$ of real analytic maps
$\psi_t:P\rightarrow\C^m$ with $\psi_0  =  \psi$ satisfying
\begin{equation*}
\begin{split}
\left(\frac{d\psi_t}{dt}\right)^b   =  (\psi_t)_*(\chi)^{a_1\ldots
\,a_{m-1}} (\Re\Omega)_{a_1 \ldots\, a_{m-1} a_m}g^{a_m b}
\end{split}
\end{equation*}

\noindent using index notation for tensors on $\C^m$. Define
$\Psi:(-\epsilon, \epsilon) \times P \rightarrow\C^m$ by $\Psi(t,p)
= \psi_t(p)$.  Then $M\! = \text{\emph{Image}}\,\Psi$ is a
nonsingular embedded (immersed) SL $m$-fold in $\C^m$.
\label{thmJoy1}
\end{thm}


In \cite[$\S 3$]{Joy3} Joyce introduces the idea of affine evolution
data with which he is able to derive an affine evolution equation,
and therefore reduces the infinite-dimensional problem of Theorem
\ref{thmJoy1} to a finite-dimensional one.

\begin{dfn}
\label{Joyaff}
Let $2\leq m \leq n$ be integers.  A set of \emph{affine evolution data} is a
pair $(P,\chi)$, where $P$ is an $(m-1)$-dimensional submanifold of $\R^n$
and $\chi:\R^n\rightarrow\Lambda^{m-1}\R^n$ is an affine map, such that
$\chi(p)$ is a nonzero element of $\Lambda^{m-1}TP$ in $\Lambda^{m-1}\R^n$ for
each nonsingular $p\in P$.  We suppose also that $P$ is not contained in any
proper affine subspace $\R^k$ of $\R^n$.

\noindent Let Aff$(\R^n,\C^m)$ be the affine space of affine maps
$\psi:\R^n\rightarrow \C^m$ and define $\mathcal{C}_P$ to be the set
of $\psi\in$Aff$(\R^n,\C^m)$ satisfying:
\begin{itemize}
\item[(i)] $\psi^*(\omega)|_P \equiv 0$,
\item[(ii)] $\psi|_{T_p P}:T_p P \rightarrow \C^m$ is injective for all $p$ in a
dense open subset of $P$.
\end{itemize}
Then (i) is a quadratic condition on $\psi$ and (ii) is an open condition on $\psi$,
 so $\mathcal{C}_P$ is a nonempty open set in the intersection of a
finite number of quadrics in Aff$(\R^n,\C^m)$.
\end{dfn}

\noindent The conditions upon $\chi$ in Definition \ref{Joyaff} are
strong.
 The result is that there are few known examples of affine evolution data.
The evolution equation derived in \cite{Joy3} is given below
\cite[Theorem 3.5]{Joy3}.

\begin{thm}\label{Joyaffthm}
Let $(P,\chi)$ be a set of affine evolution data and let
$\psi\in\mathcal{C}_P$, where $\mathcal{C}_P$ is defined in
Definition \ref{Joyaff}. Then there exist $\epsilon
> 0$ and a unique real analytic family $\{\psi_t:t\in(-\epsilon,
\epsilon)\}$ in $\mathcal{C}_P$ with $\psi_0=\psi$, satisfying
\begin{equation*}
\left(\frac{d\psi_t}{dt}(x)\right)^b = (\psi_t)_*(\chi(x))^{a_1
\ldots\, a_{m-1}}(\Re\Omega)_{a_1 \ldots\, a_{m-1} a_m} g^{a_m b}
\end{equation*}
\noindent for all $x\in\R^n$, using index notation for tensors in
$\C^m$. Furthermore $M = \{\psi_t(p):t\in(-\epsilon,\epsilon),
\hspace{2pt} p\in P\}$ is an SL $m$-fold in $\C^m$ wherever it is
nonsingular.
\end{thm}

We conclude this section by discussing the material in \cite{Joy4},
which is\linebreak particularly pertinent to $\S$\ref{second}, where
Joyce, for the majority of the paper, focuses on constructing SL
3-folds in $\C^3$ using the set of affine evolution data given below
\cite[p. 352]{Joy4}.

\begin{ex}
\label{5embed} Let $\phi: \R^2 \rightarrow \R^5$ be the embedding of
$\R^2$ in $\R^5$ given by
\begin{equation}
\begin{split}\label{phiembed}
\phi(y_1,y_2) = \left(\frac{1}{2}(y_1^2 + y_2^2), \frac{1}{2} (y_1^2 - y_2^2), y_1y_2, y_1, y_2\right).
\end{split}
\end{equation}
Then $P = \text{Image}\,\phi$ can be written as
\begin{equation*}
\begin{split}
P = \left\{ (x_1,\ldots,x_5)\in\R^5:x_1=\frac{1}{2}(x_4^2 + x_5^2), \hspace{2pt} x_2=
\frac{1}{2}(x_4^2 - x_5^2), \hspace{2pt} x_3 =x_4x_5 \right\},
\end{split}
\end{equation*}
\noindent which is diffeomorphic to $\R^2$.  From \eq{phiembed}, we
calculate, writing $e_j = \frac{\partial}{\partial x_j}\,$:
\begin{align*}
\phi_*\left(\frac{\partial}{\partial y_1}\right) & = y_1 e_1 + y_1
e_2
+ y_2 e_3 + e_4, \\[4pt]
\phi_*\left(\frac{\partial}{\partial y_2}\right) & = y_2 e_1 - y_2
e_2 + y_1 e_3 + e_5
\end{align*}
\noindent and thus
\begin{align*}
\phi_*\left(\frac{\partial}{\partial y_1} \w
\frac{\partial}{\partial y_2}\right)  =
 &\hspace{4pt}(y_1^2 + y_2^2)\hspace{1pt}e_2 \w e_3 +(y_1^2 - y_2^2)\hspace{1pt}e_1 \w e_3 - 2y_1y_2 \hspace{1pt} e_1 \w e_2
\\[-6pt]
& {} + y_1\left(e_1 \w e_5 + e_2 \w e_5 - e_3 \w e_4\right) + e_4 \w
e_5
 \\
&{} + y_2\left(-e_1 \w e_4 + e_2 \w e_4 + e_3 \w e_5 \right).
\end{align*}
Hence, if we define an affine map $\chi: \R^5 \rightarrow \Lambda^2
\R^5$ by
\begin{align}
 &\chi(x_1,\ldots,x_5) =  2x_1 e_2 \w e_3 + 2x_2 e_1 \w e_3  - 2x_3 e_1 \w e_2 + e_4 \w e_5 \nonumber\\
\label{affinchi} &{}+ x_4\left(e_1 \w e_5 + e_2 \w e_5 - e_3 \w
e_4 \right)+ x_5\left(-e_1 \w e_4 + e_2 \w e_4 + e_3 \w e_5
\right),
\end{align}
\noindent then $\chi = \phi_* \left( \frac{\partial}{\partial y_1}
\w \frac{\partial}{\partial y_2} \right)$ on $P$. Therefore $(P,
\chi)$ is a set of affine evolution data with $m=3$ and $n=5$.
\end{ex}

\noindent The main result \cite[Theorem 5.1]{Joy4} requires the
definition of a cross product $\times: \C^3 \times \C^3 \rightarrow
\C^3$, given in index notation by
\begin{equation}
\begin{split}
\label{6cross} ({\bf u} \times {\bf v})^d = ({\rm Re} \hspace{2pt}
\Omega) _{abc} {\bf u}^a {\bf v}^b g^{cd}
\end{split}
\end{equation}
\noindent for ${\bf u}, {\bf v}\in\C^3$, regarding $\C^3$ as a real
vector space.

\begin{thm}
\label{thmJoy2} Suppose that ${\bf z}_1,\ldots,{\bf
z}_6:\R\rightarrow\C^3$ are differentiable functions satisfying:
\begin{align}
\label{Joycond1}
\omega({\bf z}_2,{\bf z}_3) = \omega({\bf z}_1,{\bf z}_3) = \omega({\bf z}_1,
{\bf z}_2) &= 0, \\
\label{Joycond2}
\omega({\bf z}_1,{\bf z}_5) + \omega({\bf z}_2,{\bf z}_5) - \omega({\bf z}_3,
{\bf z}_4) &= 0, \\
\label{Joycond3}
{}-\omega({\bf z}_1,{\bf z}_4) + \omega({\bf z}_2,{\bf z}_4) + \omega({\bf z}_3,
{\bf z}_5)& = 0, \\
\label{Joycond4} \omega({\bf z}_4,{\bf z}_5) &= 0,
\end{align}
\noindent at $t=0$, and the equations:
\begin{align}
\label{Joyeq1}
\frac{d{\bf z}_1}{dt} & =  2 {\bf z}_2 \times {\bf z}_3, \\
\label{Joyeq2}
\frac{d{\bf z}_2}{dt} & =  2 {\bf z}_1 \times {\bf z}_3, \\
\label{Joyeq3}
\frac{d{\bf z}_3}{dt} & =  -2 {\bf z}_1 \times {\bf z}_2, \\
\label{Joyeq4} \frac{d{\bf z}_4}{dt} & =  {\bf z}_1 \times {\bf
z}_5
+ {\bf z}_2 \times {\bf z}_5 - {\bf z}_3 \times {\bf z}_4, \\
\label{Joyeq5} \frac{d{\bf z}_5}{dt} & =  -{\bf z}_1 \times {\bf
z}_4
+ {\bf z}_2 \times {\bf z}_4 + {\bf z}_3 \times {\bf z}_5, \\
\label{Joyeq6} \frac{d{\bf z}_6}{dt} & =  {\bf z}_4 \times {\bf
z}_5,
\end{align}
\noindent for all $t\in\R$, where $\times$ is defined by
\eq{6cross}.  Let $M\subseteq\C^3$ be defined by:
\begin{align*}
M  =  \bigg\{ &\frac{1}{2}(y_1^2 + y_2^2){\bf z}_1(t) + \frac{1}{2}(y_1^2 - y_2^2){\bf z}_2(t) + y_1y_2{\bf z}_3(t) \\
& {}+ y_1{\bf z}_4(t) + y_2{\bf z}_5(t) + {\bf z}_6(t): y_1, y_2, t
\in \R \bigg\}.
\end{align*}
\noindent Then $M$ is a special Lagrangian 3-fold in $\C^3$
wherever it is nonsingular.
\end{thm}

Joyce \cite{Joy4} solves \eq{Joyeq1}-\eq{Joyeq6} subject to the
conditions \eq{Joycond1}-\eq{Joycond4}, dividing the solutions into
cases based on the dimension of $\langle {\bf z}_1(t),{\bf z}_2(t),
{\bf z}_3(t) \rangle_{\R}$ for generic $t\in\R$.  We shall be
concerned with the case where dim$\langle {\bf z}_1(t),{\bf
z}_2(t),{\bf z}_3(t) \rangle_{\R} =3$, which forms the bulk of the
results of \cite{Joy4}.  The solutions in this case involve the {\it
Jacobi elliptic functions}, which we now give a brief description
of, following the material in \cite[Chapter VII] {Elliptic}.

For $k\in [0,1]$, the Jacobi elliptic functions, $\sn(u,k)$,
$\cn(u,k)$, $\dn(u,k)$, with modulus $k$ are the unique solutions to
the equations:
\begin{align*}
\left(\frac{d}{du}\sn(u,k)\right)^2 & =  (1 - \sn^2(u,k))(1 -
k^2\sn^2(u,k)),
\\
\left(\frac{d}{du}\cn(u,k)\right)^2 & =  (1 - \cn^2(u,k))(1 - k^2
+ k^2
\cn^2(u,k)), \\
\left(\frac{d}{du}\dn(u,k)\right)^2 & =  -(1 - \dn^2(u,k))(1 - k^2
- \dn^2(u,k)),
\end{align*}
\noindent with the initial conditions \vspace{8pt}

$\begin{array}{rrr}
\sn(0,k)=0, & \hspace{20pt} \cn(0,k)=1, & \hspace{20pt} \dn(0,k)=1, \\
\frac{d}{du}\sn(0,k)=1, & \hspace{20pt} \frac{d}{du}\cn(0,k)=0, & \hspace{20pt} \frac{d}{du}\dn(0,k)=0.
\end{array}$

\vspace{8pt}
\noindent They also satisfy the following identities and differential equations:
\begin{align*}
\sn^2(u,k)+\cn^2(u,k)&=1,\\
k^2\sn^2(u,k)+\dn^2(u,k)&=1,\\
\frac{d}{du}\sn(u,k)&=\cn(u,k)\dn(u,k),\\
\frac{d}{du}\cn(u,k)&=-\sn(u,k)\dn(u,k),\\
\frac{d}{du}\dn(u,k)&=-k^2\sn(u,k)\cn(u,k).
\end{align*}
\noindent For $k=0,1$ they reduce to familiar functions:

\vspace{8pt}
$\begin{array}{lll}
\sn(u,0)=\sin u, & \hspace{20pt} \cn(u,0)=\cos u, & \hspace{20pt}\dn(u,0)=1, \\
\sn(u,1)=\tanh u, & \hspace{20pt}\cn(u,1)=\sech u, & \hspace{20pt}\dn(u,1)=\sech u.
\end{array}$
\vspace{8pt}

\noindent For each $k\in [0,1)$ they are periodic functions.

\medskip

The embedding given in Example \ref{5embed} was constructed by
considering the \ action of $\SL(2,\R)\ltimes\R^2$ on $\R^2$. Hence,
Joyce \cite[Proposition 9.1]{Joy4} shows that solutions of
\eq{Joyeq1}-\eq{Joyeq3}, satisfying the condition \eq{Joycond1}, are
equivalent under the natural actions of $\SL(2,\R)$ and $\SU(3)$ to
a solution of the form ${\bf z}_1=(z_1,0,0)$, ${\bf z}_2=(0,z_2,0)$,
${\bf z}_3=(0,0,z_3)$, for differentiable functions
$z_1,z_2,z_3:\R\rightarrow\C$.  Therefore we assume that the
solution is of this form.  Equations \eq{Joyeq1}-\eq{Joyeq3} become:
\begin{equation}
\label{simplified}
\frac{dz_1}{dt} = 2\hspace{2pt}\overline{z_2z_3}, \hspace{20pt} \frac{dz_2}{dt}= -2\hspace{2pt}\overline{z_3z_1},
\hspace{20pt} \frac{dz_3}{dt}=-2\hspace{2pt}\overline{z_1z_2}.
\end{equation}
\noindent The next result is taken from \cite[Proposition
9.2]{Joy4}.

\begin{prop}
\label{Joyprop} Given any initial data $z_1(0),z_2(0),z_3(0)$,
solutions to \eq{simplified} exist for all $t\in\R$. Wherever the
$z_j(t)$ are nonzero they may be written as:

\vspace{8pt}
$\begin{array}{rrr}
2z_1=e^{i\theta_1}\sqrt{\alpha_1^2+v}, & \hspace{8pt} 2z_2=e^{i\theta_2}\sqrt{\alpha_2^2-v}, &
\hspace{8pt} 2z_3=e^{i\theta_3}\sqrt{\alpha_3^2-v},
\end{array}$
\vspace{8pt}

\noindent where $\alpha_j\in\R$ for all $j$ and
$v,\theta_1,\theta_2,\theta_3:\R \rightarrow\R$ are differentiable
functions.  Let $\theta=\theta_1+\theta_2+ \theta_3$ and let
$Q(v)=(\alpha_1^2+v)(\alpha_2^2-v)(\alpha_3^2-v)$.  Then there
exists $A\in\R$ such that $Q(v)^{\frac{1}{2}}\sin\theta=A$.
\end{prop}

We state the main theorem that we shall require in $\S$\ref{second},
\cite[Theorem 9.3]{Joy4}.

\begin{thm}
\label{Joythm}
Using the notation of Proposition \ref{Joyprop}, let $\alpha_j>0$ for all $j$ and
$\alpha_1^{-2} = \alpha_2^{-2} + \alpha_3^{-2}$.  Suppose that
$v$ has a minimum at $t=0$, that $\theta_2(0)=\theta_3(0)=0$, $A\geq 0$ and
that $\alpha_2\leq\alpha_3$.  Then exactly one of the following four cases holds:
\begin{rlist}
\item $A=0$ and $\alpha_2=\alpha_3$, and $z_1,z_2,z_3$ are given by:
\begin{align*}
2z_1(t) & =  \sqrt{3}\alpha_1\tanh\left(\sqrt{3}\alpha_1 t\right), \\
2z_2(t) & =  2z_3(t) = \sqrt{3}\alpha_1\sech\left(\sqrt{3}\alpha_1
t\right);
\end{align*}
\item $A=0$ and $\alpha_2<\alpha_3$, and $z_1,z_2,z_3$ are given
by:
\begin{align*}
2z_1(t) & =  \sqrt{\alpha_1^2+\alpha_2^2}\sn(\sigma t, \tau), \\
2z_2(t) & =  \sqrt{\alpha_1^2+\alpha_2^2}\cn(\sigma t, \tau), \\
2z_3(t) & =  \sqrt{\alpha_1^2+\alpha_3^2}\dn(\sigma t, \tau),
\end{align*}
\noindent where $\sigma = \sqrt{\alpha_1^2+\alpha_3^2}$ and $\tau
= \sqrt{ \frac{\alpha_1^2+\alpha_2^2}{\alpha_1^2+\alpha_3^2}};$
\vspace{8pt} \item $0<A<\alpha_1\alpha_2\alpha_3$.  Let the roots
of $Q(v)-A^2$ be $\gamma_1,\gamma_2,\gamma_3$, ordered such that
$\gamma_1\leq0\leq\gamma_2\leq \gamma_3$.  Then
$v,\theta_1,\theta_2,\theta_3$ are given by:
\begin{align*}
v(t)&=\gamma_1 + (\gamma_2-\gamma_1)\sn^2(\sigma t,\tau), \\
\theta_1(t)&=\theta_1(0) - A \int_0^t
\frac{ds}{\alpha_1^2+\gamma_1+(\gamma_2 -
 \gamma_1)\sn^2(\sigma s, \tau)}, \\
\theta_2(t)&=A \int_0^t \frac{ds}{\alpha_2^2-\gamma_1-(\gamma_2 -
 \gamma_1)\sn^2(\sigma s, \tau)}, \\
\theta_3(t)&=A \int_0^t \frac{ds}{\alpha_3^2-\gamma_1-(\gamma_2 -
 \gamma_1)\sn^2(\sigma s, \tau)},
\end{align*}
\noindent where $\sigma = \sqrt{\gamma_3-\gamma_1}$ and
$\tau=\sqrt{\frac{ \gamma_2-\gamma_1}{\gamma_3-\gamma_1}};$
\vspace{8pt} \item $A=\alpha_1\alpha_2\alpha_3$.  Define
$a_1,a_2,a_3\in\R$ by:

\begin{center}
$\begin{array}{rrr}
a_1=-\displaystyle\frac{\alpha_2\alpha_3}{\alpha_1},&\hspace{20pt}
a_2=\displaystyle\frac{\alpha_3\alpha_1}{\alpha_2},&\hspace{20pt}
a_3=\displaystyle\frac{\alpha_1\alpha_2}{\alpha_3},
\end{array}$
\end{center}

\noindent then $a_1+a_2+a_3=0$ since $\alpha_1^{-2}=\alpha_2^{-2}+
\alpha_3^{-2}$ and $z_1,z_2,z_3$ are given by:

\begin{center}
$\begin{array}{rrr} 2z_1(t)=i\alpha_1e^{ia_1t},\quad &
2z_2(t)=\alpha_2e^{ia_2t},\quad & 2z_3(t)=\alpha_3e^{ia_3t}.
\end{array}$
\end{center}
\end{rlist}
\end{thm}

\section{The First Evolution Equation}
\label{firsteq}
\setcounter{thm}{0}

To derive our evolution equation we shall require two results
related to \emph{real analyticity}. The first follows from the
minimality of associative 3-folds, as discussed in \cite{HarLaw}.

\begin{thm}
Let $N$ be an associative 3-fold in $\R^7$.  Then $N$ is real analytic wherever it is
nonsingular.
\label{thm1}
\end{thm}

The proof of the next result \cite[Theorem IV.4.1]{HarLaw} relies on
the \emph{Cartan--K\"ahler Theorem}, which is only applicable in the
real analytic category.

\begin{thm}
Let\/ $P$ be a 2-dimensional real analytic submanifold of\linebreak
$\Im\O \cong \R^7$.  Then there exists a unique real analytic
associative 3-fold $N$ in $\R^7$ which contains $P$. \label{thm2}
\end{thm}

We now formulate an evolution equation for associative 3-folds,
given a 2-dimensional real analytic submanifold of $\R^7$, following
Theorem \ref{thmJoy1}.

\begin{thm}
Let $P$ be a compact, orientable, 2-dimensional, real analytic
manifold, let $\chi$ be a real analytic nowhere vanishing section of
$\Lambda^2TP$, and let $\psi:P \rightarrow\R^7$ be a real analytic
embedding (immersion).  Then there exist $\epsilon > 0$ and a unique
family $\{\psi_t:t\in(-\epsilon,\epsilon)\}$ of real analytic maps
$\psi_t:P\rightarrow\R^7$ with $\psi_0  =  \psi$ satisfying
\begin{equation}
\begin{split}
\left(\frac{d\psi_t}{dt}\right)^d   =
(\psi_t)_*(\chi)^{ab}\varphi_{abc}g^{cd},
\end{split}
\label{evolve1}
\end{equation}

\noindent where $g^{cd}$ is the inverse of the Euclidean metric on
$\R^7$, using index notation for tensors on $\R^7$. Define
$\Psi:(-\epsilon, \epsilon) \times P \rightarrow\R^7$ by $\Psi(t,p)
= \psi_t(p)$. Then $M =$ \emph{Image} $\Psi$ is a nonsingular
embedded (immersed) associative 3-fold in $\R^7$. \label{thm3}
\end{thm}

\noindent Note that we are realising $M$ as the total space of a one
parameter family of two-dimensional manifolds $\{ P_t :
t\in(-\epsilon,\epsilon)\}$, where each $P_t$ is diffeomorphic to
$P$, satisfying a first-order ordinary differential equation in $t$
with initial condition $P_0 = P$.

\begin{proof}
Equation \eq{evolve1} is an evolution equation for maps $\psi_t:P
\rightarrow\R^7$ with the initial condition $\psi_0 = \psi$. Since
$P$ is compact and $P$, $\chi$, $\psi$ are real analytic, the {\it
Cauchy--Kowalevsky Theorem} \cite[p. 234]{Racke} from the theory of
partial differential equations gives $\epsilon
> 0$ such that a unique solution to the evolution equation
exists for $t\in(-\epsilon,\epsilon)$.

By Theorem \ref{thm2}, there exists a unique real analytic
associative 3-fold $N \subseteq \R^7$ such that $\psi(P)\subseteq
N$.  Consider a family $\{\tilde{\psi_t}:t\in(-\tilde{\epsilon},
\tilde{\epsilon})\}$, for some $\tilde{\epsilon} > 0$, of real
analytic maps $\tilde{\psi_t}:P \rightarrow N$, with $\tilde{\psi_0}
= \psi$, satisfying
\begin{equation}
\begin{split}
\left(\frac{d\tilde{\psi_t}}{dt}\right)^d   =
(\tilde{\psi_t})_*(\chi)^{ab}(\varphi | _N)_{abc} (g | _N)^{cd},
\label{evolve2}
\end{split}
\end{equation}

\noindent using index notation for tensors on $N$.  By the same
argument as above, a unique solution exists to \eq{evolve2} for some
$\tilde{\epsilon} > 0$.

Let $p \in P$, $t\in(-\tilde{\epsilon},\tilde{\epsilon})$ and set $x
= \tilde{\psi_t}(p) \in N$.  Let $V = (T_xN)^{\bot}$ in $\R^7$, so
$\R^7 = T_xN \oplus V$ and $(\R^7)^{*} = T_{x}^{*}N \oplus V^{*}$.
This induces a splitting:
\begin{equation*}
\begin{split}
\Lambda^3(\R^7)^* = \sum_{k=0}^{3} \Lambda^k T_{x}^*N \otimes \Lambda^{3-k}V^*.
\end{split}
\end{equation*}
Note that $\varphi\in\Lambda^3(\R^7)^*$ and that $N$ is calibrated
with respect to $\varphi$ as $N$ is an associative 3-fold.
Therefore, the component of $\varphi$ in $\Lambda^2 T_{x}^*N \otimes
V^*$ is zero since this measures the change in $\varphi | _{T_xN}$
under small variations of $T_xN$, but $\varphi | _{T_xN}$ is maximum
and therefore stationary.  Since $(\tilde{\psi_t})_*(\chi) | _p$
lies in $\Lambda^2 T_{x}N$, $(\tilde{\psi_t})_*(\chi)^{ab} |_p
\varphi_{abc}$ lies in $T_{x}^*N$, because the component in $V^*$
comes from the\ component of $\varphi$ in $\Lambda^2 T_{x}^*N
\otimes V^*$, which is zero by above.  Therefore,
\begin{equation*}
\begin{split}
(\tilde{\psi_t})_*(\chi)^{ab} |_p \varphi_{abc} =
(\tilde{\psi_t})_*(\chi)^{ab} |_p (\varphi | _{T_xN})_{abc}.
\end{split}
\end{equation*}
As $(\R^7)^* = T_{x}^{*}N \oplus V^*$ is an orthogonal
decomposition, $g^{cd} = (g | _{T_xN})^{cd} + h^{cd}$ for some $h\in
S^2V$.  Then $(\tilde{\psi_t})_*(\chi)^{ab} |_p (\varphi |
_{T_xN})_{abc} h^{cd}$ is zero because
$(\tilde{\psi_t})_*(\chi)^{ab} |_p (\varphi | _{T_xN})_{abc} \in
T_{x}^{*}N$ and $h\in S^2 V$, so their contraction is zero.  Hence,
\begin{equation*}
\begin{split}
(\tilde{\psi_t})_*(\chi)^{ab} \varphi_{abc}g^{cd} =
(\tilde{\psi_t})_*(\chi)^{ab}(\varphi | _N)_{abc}(g | _N)^{cd}
\end{split}
\end{equation*}
\noindent for all $p\in P$ and $t \in (-\tilde{\epsilon},
\tilde{\epsilon})$. Thus the family
$\{\tilde{\psi_t}:t\in(-\tilde{\epsilon}, \tilde{\epsilon})\}$
satisfies \eq{evolve1} and $\tilde{\psi_0} = \psi$, which implies
that $\tilde{\psi_t} = \psi_t$ by uniqueness.

Hence, $\psi_t$ maps $P$ to $N$ and $\Psi$ maps $(-\epsilon,
\epsilon) \times P$ to $N$ for $\epsilon$ sufficiently small.
Suppose $\psi$ is an embedding. Then $\psi_t:P \rightarrow N$ is an
embedding for small $t$. Moreover, $\frac{d\psi_t}{dt}$ is a normal
vector field to $\psi_t(P)$ in $N$ with length $| (\psi_t)_*(\chi)
|$, so, since $\chi$ is nowhere vanishing, this vector field is
nonzero. We deduce that $\Psi$ is an embedding for small $\epsilon$,
with Image $\Psi = M$ an open subset of $N$, and conclude that $M$
is an associative 3-fold.
 Similarly if $\psi$ is an immersion.
\end{proof}

\section{The Second Evolution Equation}
\setcounter{thm}{0}
\label{second}

In general it is difficult to use Theorem \ref{thm3} as stated to
construct associative 3-folds, since it is an {\it
infinite-dimensional} evolution problem.  We follow the material in
\cite[$\S 3$]{Joy3} to reduce the theorem to a {\it
finite-dimensional} problem.

\begin{dfn}
Let $n \geq 3$ be an integer.  A set of \emph{affine evolution data}
is a pair $(P, \chi)$, where $P$ is a 2-dimensional submanifold of
$\R^n$ and $\chi: \R^n\rightarrow \Lambda^2\R^n$ is an affine map,
such that $\chi(p)$ is a nonzero element of $\Lambda^2TP$ in
$\Lambda^2\R^n$ for each nonsingular point $p\in P$.  Further,
suppose that $P$ is not contained in any proper affine subspace
$\R^k$ of $\R^n$.

\noindent Let Aff$(\R^n, \R^7)$ be the affine space of affine maps
$\psi:\R^n\rightarrow \R^7$.  Define $\mathcal{C}_P$ as the set of
$\psi\in\text{Aff}(\R^n, \R^7)$ such that $\psi|_{T_p P} : T_p P
\rightarrow \R^7$ is injective for all $p$ in a dense open subset of
$P$.  Let $M$ be an associative 3-plane in $\R^7$. Then generic
linear maps $\psi:\R^n\rightarrow M$ will satisfy the condition to
be members of $\mathcal{C}_P$. Hence $\mathcal{C}_P$ is non-empty.
\label{affine}
\end{dfn}

\noindent We formulate our second evolution equation following
Theorem \ref{Joyaffthm}.

\begin{thm}
\label{thm4} Let $(P, \chi)$ be a set of affine evolution data and
$n$, \emph{Aff}$(\R^n,\R^7)$ and $\mathcal{C}_P$ be as in Definition
\ref{affine}.  Suppose $\psi \in \mathcal{C}_P$.  Then there exist
$\epsilon > 0$ and a unique one parameter family
$\{\psi_t:t\in(-\epsilon,\epsilon)\} \subseteq \mathcal{C}_P$ of
real analytic maps with $\psi_0 = \psi$ satisfying
\begin{equation}
\begin{split}
\label{evolvefin}
\left(\frac{d\psi_t}{dt}(x)\right)^d = (\psi_t)_* (\chi(x))^{ab} \varphi_{abc}
g^{cd}
\end{split}
\end{equation}
\noindent for all $x\in\R^n$, using index notation for tensors on $\R^7$, where
$g^{cd}$ is the inverse of the Euclidean metric on $\R^7$.  Define $\Psi:(-
\epsilon,\epsilon) \times P \rightarrow \R^7$ by $\Psi(t,p) = \psi_t(p)$.  Then
$M =$ \emph{Image} $\Psi$ is an associative 3-fold wherever it is nonsingular.
\end{thm}

\begin{proof}
It is sufficient to restrict to the case of linear maps $\psi:\R^n
\rightarrow \R^7$ since $\R^n$ can be regarded as $\R^n\times \{1\}
\subseteq \R^{n+1} = \R^n \times \R$, and therefore any affine map\
$\psi:\R^n\rightarrow\R^7$ can be uniquely extended to a linear map
$\tilde{\psi}:\R^{n+1} \rightarrow \R^7$. We\ denote the space of
linear maps from $\R^n$ to $\R^7$ by Hom$(\R^n, \R^7)$.  Therefore\
\eq{evolvefin} is a well-defined first-order ordinary differential
equation upon the maps $\psi_t \in$ Hom$(\R^n,\R^7)$ of the form
$\frac{d\psi_t}{dt} = Q(\psi_t)$, where $Q$ is a quadratic.  Hence,
by the theory of ordinary differential equations, there exist
$\epsilon >0$ and a unique real analytic family
$\{\psi_t:t\in(-\epsilon,\epsilon)\} \subseteq$ Hom$(\R^n, \R^7)$,
with $\psi_0 = \psi$, satisfying equation \eq{evolvefin}.

Having established existence and uniqueness we can then follow the
proof of Theorem \ref{thm3}, noticing that we may drop the
assumption made there of the compactness of $P$, since it was only
used to establish the existence of the required family of maps. Note
that \eq{evolve1} is precisely the restriction of \eq{evolvefin} to
$x\in P$, so we deduce that $M$ is an associative 3-fold wherever it
is nonsingular.

We need only show now that the family constructed lies in
$\mathcal{C}_P$. Note that the requirement that $\psi_t|_{T_p P}:
T_p P \rightarrow \R^7$ is injective for all $p$ in an open dense
subset of $P$ is clearly an open condition, and that it holds at
$\psi_0 = \psi$ since $\psi \in \mathcal{C}_P$. Thus, by selecting a
sufficiently small value of $\epsilon$, we see that $\psi_t \in
\mathcal{C}_P$ for all $t \in (-\epsilon, \epsilon)$ and the proof
is complete.
\end{proof}

Before we construct associative 3-folds using this result, it is
worth noting that using quadrics to provide affine evolution data as
in \cite{Joy3} would not be a worthwhile enterprise. Suppose
$Q\subseteq\R^3$ is a quadric and that $L:\R^3\rightarrow\R^7$ is a
linear map. Then we can transform $\R^7$ using $\text{G}_2$ such
that, if we write $\R^7 = \R \oplus \C^3$, then $L(\R^3) \subseteq
\C^3$ is a Lagrangian plane. Therefore, evolving $Q$ using
\eq{evolvefin} will only produce SL 3-folds, which have already been
studied in \cite{Joy3}.

\medskip

Let us now return to the affine evolution data given in Example
\ref{5embed} and use Theorem \ref{thm4} to construct associative
3-folds. Let $(P,\chi)$ be as in Example \ref{5embed} and define
affine maps $\psi_t:\R^5\rightarrow\R^7$ by:
\begin{equation}
\begin{split}
\label{psiaff}
\psi_t(x_1,\ldots,x_5) = {\bf w}_1(t)x_1 + \ldots + {\bf w}_5(t)x_5 + {\bf w}_6 (t),
\end{split}
\end{equation}
\noindent where ${\bf w}_j:\R \rightarrow \R^7$ are smooth functions
for all $j$.  Using the notation of Example \ref{5embed}, we see
that $(\psi_t)_*(e_j) = {\bf w}_j$ for $j=1,\ldots,5$.  Hence, by
equation \eq{affinchi} for $\chi$, equation \eq{cross2} for the
cross product on $\R^7$ and \eq{evolvefin} we have that
\begin{align}
\frac{d\psi_t}{dt}(x_1,\ldots,x_5)& = 2x_1 {\bf w}_2 \times {\bf
w}_3 + 2x_2 {\bf w}_1 \times {\bf w}_3 - 2 x_3 {\bf w}_1 \times
{\bf w}_2
\nonumber\\
& + x_4({\bf w}_1 \times {\bf w}_5 + {\bf w}_2 \times {\bf w}_5 -
{\bf w}_3 \times {\bf w}_4)\nonumber\\[2pt]
\label{dpsiaff} & + x_5(-{\bf w}_1 \times {\bf w}_4 + {\bf w}_2
\times {\bf w}_4 + {\bf w}_3 \times {\bf w}_5)+{\bf w}_4\times{\bf
w}_5
\end{align}
\noindent for all $(x_1,\ldots,x_5) \in \R^5$.  Therefore, from
\eq{psiaff} and
 \eq{dpsiaff} we get the following result.

\begin{thm}
\label{evolvethm} Let ${\bf w}_1,\ldots,{\bf w}_6:\R\rightarrow\R^7$
be differentiable functions satisfying
\begin{align}
\label{ch5.1}
\frac{d{\bf w}_1}{dt} & =  2 {\bf w}_2 \times {\bf w}_3, \\
\label{ch5.2}
\frac{d{\bf w}_2}{dt} & =  2 {\bf w}_1 \times {\bf w}_3, \\
\label{ch5.3}
\frac{d{\bf w}_3}{dt} & =  -2 {\bf w}_1 \times {\bf w}_2, \\
\label{ch5.4}
\frac{d{\bf w}_4}{dt} & =   {\bf w}_1 \times {\bf w}_5 + {\bf w}_2 \times {\bf w}_5 - {\bf w}_3 \times {\bf w}_4,  \\
\label{ch5.5}
\frac{d{\bf w}_5}{dt} & =  - {\bf w}_1 \times {\bf w}_4 + {\bf w}_2 \times {\bf w}_4 + {\bf w}_3 \times {\bf w}_5,  \\
\label{ch5.6} \frac{d{\bf w}_6}{dt} & =   {\bf w}_4 \times {\bf
w}_5.
\end{align}
\noindent Then $M$, given by:
\begin{align*}
M  =   \bigg\{ &\frac{1}{2}(y_1^2 + y_2^2){\bf w}_1(t) + \frac{1}{2}(y_1^2 - y_2^2){\bf w}_2(t) + y_1y_2{\bf w}_3(t) \\
&+ y_1{\bf w}_4(t) + y_2{\bf w}_5(t) + {\bf w}_6(t): y_1, y_2,t \in
\R \hspace{2pt} \bigg\},
\end{align*}
\noindent is an associative 3-fold in $\R^7$ wherever it is
nonsingular.
\end{thm}

\noindent Theorem \ref{thm4} only gives us that the associative
3-fold $M$ is defined for $t$ in some small open neighbourhood of
zero, but work later in this section shows that $M$ is indeed
defined for all $t$ as stated in the above theorem.

The equations we have just obtained fall naturally into three parts:
 \eq{ch5.1}-\eq{ch5.3} show that ${\bf w}_1, {\bf w}_2, {\bf w}_3$
 evolve amongst themselves; \eq{ch5.4}-\eq{ch5.5} are {\it linear}
 equations for ${\bf w}_4$ and ${\bf w}_5$ once ${\bf w}_1,{\bf w}_2,{\bf w}_3$
are known; and \eq{ch5.6} defines ${\bf w}_6$ once the functions
${\bf w}_4$ and ${\bf w}_5$ are known.  Moreover, these equations
are very similar to \eq{Joyeq1}-\eq{Joyeq6}, given in Theorem
\ref{thmJoy2}, the only difference being that here our functions and
cross products are defined on $\R^7$ rather than $\C^3$. If we could
show that any solutions ${\bf w}_1, {\bf w}_2, {\bf w}_3$ are
equivalent to functions ${\bf z}_1, {\bf z}_2, {\bf z}_3$, lying in
$\C^3$, satisfying \eq{Joyeq1}-\eq{Joyeq3} and \eq{Joycond1}, then
we would be able to use results from \cite{Joy4} to hopefully
construct associative 3-folds which are not SL 3-folds.  It is to
this end that we now proceed.

\medskip

Suppose that ${\bf w}_1(t), {\bf w}_2(t), {\bf w}_3(t)$ are
solutions to \eq{ch5.1}-\eq{ch5.3}. Let $w_j = {\bf w}_j(0)$ for all
$j$ and let $v = [w_1, w_2, w_3]$, as defined by \eq{asctr}.

If $v=0$, then, by Proposition \ref{ass}, $\langle w_1,w_2,w_3
\rangle_{\R}$ lies in an associative 3-plane which we can map to
$\R^3 \subseteq \C^3 \subseteq \R^7 = \R \oplus \C^3$, since
$\text{G}_2$ acts transitively on associative 3-planes \cite[Theorem
IV.1.8]{HarLaw}. Let $z_1,z_2,z_3$ be the images of $w_1,w_2,w_3$
under this transformation and let $\omega$ be the standard
symplectic form on $\C^3$, which in terms of coordinates
$(x_1,\ldots,x_7)$ on $\R^7$ is given by:
\begin{equation*}
\begin{split}
\omega = dx_2 \w dx_3 + dx_4 \w dx_5 + dx_6 \w dx_7.
\end{split}
\end{equation*}
\noindent Then, $z_1,z_2,z_3$ lie in $\R^3 \subseteq \C^3$ and so
$\omega(z_j,z_k)=0$ for $j \neq k$.

If $v\neq 0$, then $v$ is orthogonal to $w_j$ for all $j$ by
Proposition \ref{assprops}, so we can split $\R^7 = \R \oplus \C^3$
where $\R = \langle v \rangle$ and $\C^3 = {\langle v
\rangle}^{\perp}$.  Hence, $w_j$ lies in $\C^3$ for all $j$ with
respect to this splitting.  By Proposition \ref{assprops}, $v$ is
orthogonal to $[w_j,w_k] = w_jw_k - w_kw_j = 2w_j\times w_k$ and
therefore, from \eq{innerphi},
\begin{equation*}
\begin{split}
\varphi_{abc}v^a {w_j}^b {w_k}^c = 0
\end{split}
\end{equation*}
\noindent using index notation for tensors on $\R^7$.  Note that we
can write:
\begin{equation}
\begin{split}
\label{splitting}
\varphi = dx_1\w\omega + {\rm Re} \hspace{2pt} \Omega,
\end{split}
\end{equation}
\noindent where $\Omega$ is the holomorphic volume form on $\C^3$.
Therefore, $\varphi_{abc} v^a = |v| \omega_{bc}$ and hence, since
$|v| \neq 0$, $\omega(w_j,w_k)=0$.

From equations \eq{cross2} and \eq{6cross} defining the cross
products on $\R^7$ and $\C^3$ respectively and \eq{splitting} above,
we see that, for vectors ${\bf x}, {\bf y}\in\C^3\subseteq\R^7$,
\begin{equation}
\begin{split}
\label{cross67}
{\bf x}\times{\bf y} = {\bf x}\times^{\prime}{\bf y} +
\omega({\bf x},{\bf y}){\bf e}_1,
\end{split}
\end{equation}
\noindent where $\times^{\prime}$ is the cross product on $\C^3$ and
${\bf e}_1= (1,{\bf 0})\in\R\oplus\C^3=\R^7$.  We have shown that,
using a $\text{G}_2$ transformation, we can map the solutions ${\bf
w}_1(t), {\bf w}_2(t), {\bf w}_3(t)$ to solutions ${\bf z}_1(t),
{\bf z}_2(t), {\bf z}_3(t)$ such that ${\bf z}_j(0)\in\C^3
\subseteq\R^7$ and $\omega({\bf z}_j(0),{\bf z}_k(0)) = 0$.  Our
remarks above about \eq{ch5.1}-\eq{ch5.3}, and the relationship
\eq{cross67} between the cross products on $\C^3$ and $\R^7$, show
that ${\bf z}_1(t), {\bf z}_2(t), {\bf z}_3(t)$ must remain in
$\C^3$ and satisfy \eq{Joyeq1}-\eq{Joyeq3} along with\linebreak
condition \eq{Joycond1}.  Hence, any solution of
\eq{ch5.1}-\eq{ch5.3} is equivalent up to a $\text{G}_2$
transformation to a solution to the corresponding equations in
Theorem \ref{thmJoy2}.

\medskip

We now perform a parameter count in order to calculate the dimension
of the family of associative 3-folds constructed by Theorem
\ref{evolvethm}. The initial data ${\bf w}_1(0),\ldots,{\bf w}_6(0)$
has 42 real parameters, which implies that dim $\mathcal{C}_P=42$
(using the notation of Definition \ref{affine}), and so the family
of curves in $\mathcal{C}_P$ has dimension 41, which corresponds to
factoring out translation in $t$.
  It is shown in \cite{Joy4} that $\GL(2,\R)\ltimes\R^2$ acts on this family of curves and, because of the internal symmetry of the evolution data, any two
curves related by this group action give the same 3-fold.  Therefore we have
to reduce the dimension of distinct associative 3-folds up to this group action
 by 6 to 35.  We can also identify any two associative 3-folds which are
isomorphic under automorphisms of $\R^7$, i.e. up to the action of
$\text{G}_2 \ltimes \R^7$, and so we reduce the dimension by 21 to
14.

In conclusion, the family of associative 3-folds constructed in this
section has dimension 14, whereas the dimension of the family of SL
3-folds constructed in Theorem \ref{thmJoy2} has dimension 9, so not
only do we know that we have constructed new geometric objects, but
also how many more interesting parameters we expect to find.

\subsection{Singularities of these associative 3-folds}
We study the singularities of the 3-folds constructed by Theorem
\ref{evolvethm} by introducing the function $F:\R^3\rightarrow\R^7$
defined by:
\begin{align}
F(y_1,y_2,t)  =&  \textstyle\frac{1}{2}\displaystyle (y_1^2 +
y_2^2){\bf w}_1(t) + \textstyle\frac{1}{2}\displaystyle(y_1^2 -
y_2^2){\bf w}_2(t)
 + y_1y_2{\bf w}_3(t) \nonumber \\
\label{singmap} &  {}+ y_1{\bf w}_4(t) + y_2{\bf w}_5(t) + {\bf
w}_6(t).
\end{align}
\noindent  Clearly, $F$ is smooth and, if
$dF|_{(y_1,y_2,t)}:\R^3\rightarrow\R^7$ is injective for all
$(y_1,y_2,t)\in\R$, then $F$ is an immersion and $M
=\text{Image}\,F$ is nonsingular.  Therefore the possible
singularities of $M$ correspond to points where $dF$ is not
injective.
  Since we have from \eq{ch5.1}-\eq{ch5.6} that
\begin{equation*}
\begin{split}
\frac{\partial F}{\partial y_1} \times \frac{\partial F}{\partial y_2} =
\frac{\partial F}{\partial t},
\end{split}
\end{equation*}
\vspace{-16pt}

\noindent $\frac{\partial F}{\partial t}$ is perpendicular to the
other two partial derivatives and it is zero if and only if the
$y_1$ and $y_2$ partial derivatives are linearly dependent. We
deduce that $F$ is an immersion if and only if $\frac{\partial
F}{\partial y_1}$ and $\frac{\partial F}{\partial y_2}$ are linearly
 independent, since $dF$ is injective if and only if the three partial
derivatives of $F$ are linearly independent.  The condition for $F$
to be an immersion at $(0,0,0)$ is that ${\bf w}_4(0)$ and ${\bf
w}_5(0)$ are linearly independent.

\medskip

We perform a parameter count for the family of singular associative
3-folds constructed by Theorem \ref{evolvethm}.  The set of initial
data ${\bf w}_1(0),\ldots,{\bf w}_6(0)$, with ${\bf w}_4(0)$ and
${\bf w}_5(0)$ linearly dependent, has dimension $28+8=36$, since
the set of linearly dependent pairs in $\R^7$ has dimension 8.  We
saw in the earlier parameter count above that the set of initial
data without any restrictions had dimension 42. Hence, the condition
that $F$ is not an immersion at $(0,0,0)$ is of real codimension 6,
but this is clearly true for any point in $\R^3$ and therefore it is
expected that the family of singular associative 3-folds will be of
codimension $6-3=3$ in the family of all associative 3-folds
constructed by Theorem \ref{evolvethm}. Therefore the family of
distinct singular associative 3-folds up to automorphisms of $\R^7$
should have dimension $14-3=11$.  Thus generic associative 3-folds
constructed by Theorem \ref{evolvethm} will be nonsingular.
Moreover, the dimension of the family of singular associative
3-folds is greater than the dimension of the family of singular SL
3-folds constructed from the same evolution data (which has
dimension 8).

\medskip

We now model $M=$ Image $F$ near a singular point, which we take to
be the origin without loss of generality.  Therefore, we expand
${\bf w}_1(t),\ldots,{\bf w}_6(t)$ about $t=0$ to study the
singularity. Since $dF$ is not injective at the origin, ${\bf
w}_4(0)$ and ${\bf w}_5(0)$ are linearly dependent.  As mentioned
above, Joyce \cite[$\S$5.1]{Joy4} describes how internal symmetry of
the evolution data gives rise to an action of $\GL(2,\R)\ltimes\R^2$
upon ${\bf w}_1(t),\ldots,{\bf w}_6(t)$, under which the associative
3-fold constructed is invariant.  A rotation of $\R^2$ by an angle
$\theta$ transforms ${\bf w}_4(0)$ and ${\bf w}_5(0)$ to
\begin{align*}
{\bf \tilde{w}}_4(0) & =  \cos\theta {\bf w}_4(0) - \sin\theta {\bf w}_5(0), \\
{\bf \tilde{w}}_5(0) & =  \sin\theta {\bf w}_4(0) + \cos\theta {\bf
w}_5(0).
\end{align*}
\noindent Since ${\bf w}_4(0)$ and ${\bf w}_5(0)$ are linearly
dependent, $\theta$ may be chosen so that ${\bf \tilde{w}}_5(0)=0$.
We may therefore suppose that ${\bf w}_5(0)=0$ and take our initial
data to be:

\vspace{8pt}
$\begin{array}{ll}
\hspace{52pt} {\bf w}_1(0) = {\bf v} + {\bf w}, & \hspace{40pt} {\bf w}_2(0) = {\bf v} - {\bf w}, \\
\hspace{52pt} {\bf w}_3(0) = {\bf x}, &  \hspace{40pt} {\bf w}_4(0) = {\bf u}, \\
\hspace{52pt} {\bf w}_5(0) = {\bf w}_6(0) = 0,
\end{array}$
\vspace{8pt}

\noindent for vectors ${\bf u},{\bf v},{\bf w},{\bf x}\in\R^7$.
Expanding our solutions to \eq{ch5.1}-\eq{ch5.6} to low order in
$t$:
\begin{align*}
{\bf w}_1(t) &= {\bf v} + {\bf w} + 2t({\bf v}-{\bf w})\times {\bf x} + O(t^2), \\
{\bf w}_2(t) &= {\bf v} - {\bf w} + 2t({\bf v}+{\bf w})\times {\bf x} + O(t^2), \\
{\bf w}_3(t) &= {\bf x} + 4t {\bf v}\times {\bf w} + O(t^2), \\
{\bf w}_4(t) &= {\bf u} + t {\bf u}\times {\bf x} +O(t^2), \\
{\bf w}_5(t) &= 2t {\bf u}\times {\bf w} + 8t^2 {\bf x} \times
({\bf u}\times
{\bf w}) + O(t^3) \\
{\bf w}_6(t) &= 10t^3 {\bf u}\times({\bf x}\times({\bf u}\times {\bf
w})) + O(t^4).
\end{align*}
Calculating $F(y_1,y_2,t)$ near the origin, we see that the dominant
terms in the expansion are dependent upon ${\bf w}_1,{\bf w}_2,{\bf
w}_3$, which we have shown to be equivalent under $\text{G}_2$ to
solutions as given in Theorem \ref{thmJoy2}. Following Joyce
\cite[p. 363-364]{Joy4}, we consider $F(\epsilon^2 y_1,\epsilon
y_2,\epsilon t)$ for small $\epsilon$, which is given by:
\begin{align}
F(\epsilon^2y_1,\epsilon y_2,\epsilon t) = \epsilon^2 [&(y_1 +
\frac{1}{4} g({\bf u},{\bf w})t^2){\bf u} +
(y_2^2-\frac{1}{4}|{\bf u}|^2t^2){\bf w} +
2y_2t{\bf u}\times{\bf w}] &\nonumber \\
+ \epsilon^3 [ &4y_2^2t {\bf x}\times {\bf w} + y_1y_2 {\bf x} +
y_1 t {\bf u}\times {\bf x} +
8y_2 t^2 {\bf x}\times ({\bf u}\times {\bf w}) \nonumber \\
&+10t^3 {\bf u} \times ( {\bf x} \times ( {\bf u} \times {\bf
w}))] + O(\epsilon^4).\label{sing1}
\end{align}
\noindent Here we have assumed that $\omega({\bf u},{\bf w})=0$ in
order to simplify the coefficient of ${\bf u}$.  The $\epsilon^2$
terms in \eq{sing1} give us the lowest order description of the
singularity.  If we suppose that ${\bf u}$ and ${\bf w}$ are
linearly independent, which will be true in the generic case, then
${\bf u}$, ${\bf w}$ and ${\bf u}\times{\bf w}$ are linearly
independent and therefore generate an SL $\R^3$. Hence, near the
origin to lowest order, $M$ is the image of the map from $\R^3$ to
$\R^3$ given by
\begin{equation}
\begin{split}
\label{sing2}
(y_1,y_2,t)\mapsto (y_1+\frac{1}{4}g({\bf u},{\bf w})t^2,y_2^2-\frac{1}{4}|{\bf u}|^2t^2,2y_2t).
\end{split}
\end{equation}
\noindent Note that the first coordinate axis is fixed under
\eq{sing2} and, moreover, $y_2$ and $t$ are allowed to take either
sign. Therefore, \eq{sing2} is a \emph{double cover} of an SL $\R^3$
which is branched over the first coordinate axis. This is the same
behaviour as occurs in the SL case \cite[p. 364]{Joy4}.

In order to study the singularity further we consider the
$\epsilon^3$ terms in \eq{sing1}. It is generally not possible to
simplify the final cross product in the $\epsilon^3$ terms to give a
neat expression using only four vectors.  However, suppose we choose
${\bf x}=(1,\mathbf{0})\in\R\oplus\C^3\cong\R^7$ and $\{{\bf u},{\bf
w},{\bf u}\times{\bf w}\}$ to be the usual oriented orthonormal
basis for the standard $\R^3$ in $\C^3\subseteq\R^7$. Then, using
\eq{sing1} and \eq{sing2}, the next order of the singularity is the
image of the following map from $\R^3$ to $\R^7$:
\begin{equation*}
\begin{split}
\label{sing3} (y_1,y_2,t)\mapsto (\epsilon y_1y_2,\, y_1,\,-\epsilon
y_1,\, y_2^2-\frac{1}{4}\,t^2,\, 4\epsilon y_2^2 t+10\epsilon t^3,\,
2y_2t,\, 8\epsilon y_2 t^2).
\end{split}
\end{equation*}
Note that the singularity does not lie within $\C^3\subseteq\R^7$
and so we have a model for a singularity which is different from the
SL case.

\subsection{Solving the equations}

From the work above, any solution ${\bf w}_1(t),{\bf w}_2(t), {\bf
w}_3(t)$ in $\R^7$ to \eq{ch5.1}-\eq{ch5.3} is equivalent under a
$\text{G}_2$ transformation to a solution ${\bf z}_1(t),{\bf
z}_2(t),{\bf z}_3(t)$ in $\C^3$ to \eq{Joyeq1}-\eq{Joyeq3}
satisfying \eq{Joycond1}.  We can thus use results from \cite{Joy4}
to produce some associative 3-folds. However, we must exercise some
caution: we require that $\langle {\bf z}_1(t),{\bf z}_2(t),{\bf
z}_3(t):t\in\R \rangle_\R = \C^3$.  If this does not occur, there
may be a further $\text{G}_2$ transformation that preserves the
subspace spanned by the ${\bf z}_j(t)$, but transforms $\C^3$ so
that ${\bf w}_4$
 and ${\bf w}_5$ are mapped into $\C^3$, and thus the submanifold
 constructed will be an SL 3-fold embedded in $\R^7$.

When $\text{dim}\,\langle {\bf z}_1(t),{\bf z}_2(t),{\bf z}_3(t)
\rangle_\R < 3$, for generic $t\in\R$, the ${\bf z}_j(t)$ define a
subspace of an SL $\R^3$ in $\C^3$, which corresponds to an
associative $\R^3$ in $\R^7$.  The subgroup of $\text{G}_2$
preserving an associative $\R^3$ is $\SO(4)$ \cite[Theorem
IV.1.8]{HarLaw}, and the subgroup of $\SU(3)$, which is the
automorphism group of $\C^3$, preserving the standard $\R^3$ is
$\SO(3)$.  Hence, the family of different ways of identifying
$\R^7\cong\R\oplus\C^3$ such that $\langle {\bf z}_1(t),{\bf
z}_2(t),{\bf z}_3(t)\rangle_\R$ is mapped into the standard $\R^3$
in $\C^3$ contains $\SO(4) / \SO(3) \cong \mathcal{S}^3$. We
therefore have sufficient freedom left in using the $\text{G}_2$
symmetry, after mapping ${\bf w}_1,{\bf w}_2,{\bf w}_3$ into $\C^3$,
to map ${\bf w}_4$ and ${\bf w}_5$ into $\C^3$ as well. This means
that these cases will only produce SL 3-folds.

It is also true in (i) and (ii) of Theorem \ref{Joythm} that the
solutions ${\bf z}_j(t)$ define a subspace of an SL $\R^3$ in $\C^3$
and so these cases will not provide any new associative 3-folds
either. Therefore we need only consider (iii) and (iv) in Theorem
\ref{Joythm}.

\medskip

Suppose we are in the situation of Theorem \ref{Joythm} so that, if
we write $\R^7=\R\oplus\C^3$, ${\bf w}_1=(0,w_1,0,0)$, ${\bf w}_2=
(0,0,w_2,0)$, ${\bf w}_3=(0,0,0,w_3)$ for differentiable functions
$w_1,w_2, w_3:\R\rightarrow\C$.  Let ${\bf w}_4=(y,p_1,p_2,q_3)$ and
${\bf w}_5=(-x,q_1,-q_2,p_3)$, where all the various functions
defined here are differentiable.  Equations \eq{ch5.4}-\eq{ch5.5}
become
\begin{align}
\label{9.2.1}
\frac{dx}{dt} & =  {\rm Im}(\bar{w}_1p_1 - \bar{w}_2p_2 - \bar{w}_3p_3), \\
\label{9.2.2}
\frac{dp_1}{dt} & =  ixw_1 + \overline{w_2p_3} + \overline{w_3p_2}, \\
\label{9.2.3}
\frac{dp_2}{dt} & =  ixw_2 - \overline{w_3p_1} - \overline{w_1p_3}, \\
\label{9.2.4}
\frac{dp_3}{dt} & =  ixw_3 - \overline{w_1p_2} - \overline{w_2p_1}; \\
\label{9.2.5}
\frac{dy}{dt} & =  {\rm Im}(\bar{w}_1q_1 - \bar{w}_2q_2 - \bar{w}_3q_3), \\
\label{9.2.6}
\frac{dq_1}{dt} & =  iyw_1 + \overline{w_2q_3} + \overline{w_3q_2}, \\
\label{9.2.7}
\frac{dq_2}{dt} & =  iyw_2 - \overline{w_3q_1} - \overline{w_1q_3}, \\
\label{9.2.8} \frac{dq_3}{dt} & =  iyw_3 - \overline{w_1q_2} -
\overline{w_2q_1}.
\end{align}
Note that the equations on $(x,p_1,p_2,p_3)$ are the same as on
$(y,q_1,q_2,q_3)$.  Moreover, $(x,p_1,p_2,p_3)= (0,w_1,w_2,w_3)$
gives an automatic solution to \eq{9.2.1}-\eq{9.2.4} and
$(y,q_1,q_2,q_3)=(0,w_1,w_2,w_3)$ solves \eq{9.2.5}-\eq{9.2.8}.

If we write ${\bf w}_6=(z,r_1,r_2,r_3)$, where $z:\R\rightarrow\R$
and $r_1,r_2,r_3:\R\rightarrow\C$ are differentiable functions,
\eq{ch5.6} becomes
\begin{align}
\label{9.2.9}
\frac{dz}{dt} & =  {\rm Im}(\bar{p}_1q_1 - \bar{p}_2q_2 - \bar{p}_3q_3), \\
\label{9.2.10}
\frac{dr_1}{dt} & =  ixp_1 + iyq_1 +\overline{p_2p_3} + \overline{q_2q_3}, \\
\label{9.2.11}
\frac{dr_2}{dt} & =  ixp_2 - iyq_2 -\overline{p_3p_1} + \overline{q_3q_1}, \\
\label{9.2.12} \frac{dr_3}{dt} & =  ixq_3 + iyp_3
-\overline{p_1q_2} - \overline{p_2q_1}.
\end{align}
Note that the conditions that $x,y,z$ are constant correspond to
\eq{Joycond3}, \eq{Joycond2} and \eq{Joycond4} in Theorem
\ref{thmJoy2} respectively.  Calculation using \eq{9.2.1}-\eq{9.2.4}
gives
\begin{equation*}
\begin{split}
\frac{d^2x}{dt^2} = x(|w_1|^2-|w_2|^2-|w_3|^2).
\end{split}
\end{equation*}
Suppose that $x$ is a nonzero constant.  Then
$|w_1|^2-|w_2|^2-|w_3|^2\equiv 0$.  Using \eq{innerphi},
\eq{ch5.1}-\eq{ch5.3} and the alternating properties of $\varphi$:
\begin{align*}
\frac{d}{dt}\,(|w_1|^2-|w_2|^2&-|w_3|^2)=2g\!\left(\frac{d\textbf{w}_1}{dt},\textbf{w}_1\!\right)-
2g\!\left(\frac{d\textbf{w}_2}{dt},\textbf{w}_2\!\right)-2g\!\left(\frac{d\textbf{w}_3}{dt},\textbf{w}_3\!\right)\\
&=4(g(\textbf{w}_2\times\textbf{w}_3,\textbf{w}_1)-g(\textbf{w}_1\times\textbf{w}_3,\textbf{w}_2)+
g(\textbf{w}_1\times\textbf{w}_2,\textbf{w}_3))\\
&=4(\varphi(\textbf{w}_2,\textbf{w}_3,\textbf{w}_1)-\varphi(\textbf{w}_1,\textbf{w}_3,\textbf{w}_2)+
\varphi(\textbf{w}_1,\textbf{w}_2,\textbf{w}_3))\\
&=12\varphi(\textbf{w}_1,\textbf{w}_2,\textbf{w}_3).
\end{align*}
Therefore
$\varphi(\textbf{w}_1,\textbf{w}_2,\textbf{w}_3)=\Re\,(w_1w_2w_3)\equiv
0$, which occurs if and only (iv) of Theorem \ref{Joythm} holds.
However, in case (iv),
$|w_1|^2-|w_2|^2-|w_3|^2=\alpha_1^2-\alpha_2^2-\alpha_3^2$, which,
together with the condition
$\alpha_1^{-2}=\alpha_2^{-2}+\alpha_3^{-2}$, forces $\alpha_j=0$ for
all $j$ which is a contradiction.  Hence, if $x$ is constant then
$x$ has to be zero, and we have a similar result for $y$. Therefore
\eq{Joycond2}-\eq{Joycond4} correspond to $x=y=0$ and $z$ constant.
This is unsurprising since having $x=y=0$ and $z$ constant
corresponds to ${\bf w}_4,{\bf w}_5,{\bf w}_6$ remaining in $\C^3$
and thus the associative 3-fold $M$ constructed will be SL and hence
satisfy $\omega|_M\equiv0$.

\medskip

Following the discussion earlier in this subsection we consider
(iii) and (iv) of Theorem \ref{Joythm}.  However, no solutions are
known in case (iii), so we focus on case (iv).  We therefore let
$\alpha_1,\alpha_2,\alpha_3$ be positive
 real numbers satisfying $\alpha_1^{-2}=\alpha_2^{-2}+\alpha_3^{-2}$ and define $a_1,a_2,a_3$ by:
\begin{equation}
\label{as}
a_1=-\frac{\alpha_2\alpha_3}{\alpha_1}, \hspace{20pt}a_2=\frac{\alpha_3\alpha_1}{\alpha_2}, \hspace{20pt}a_3=\frac{\alpha_1\alpha_2}{\alpha_3}.
\end{equation}
\noindent By Theorem \ref{Joythm}, we have that

\vspace{8pt} $\begin{array}{rrr} 2w_1(t)=i\alpha_1e^{ia_1t}, &
\hspace{10pt} 2w_2(t)=\alpha_2e^{ia_2t}, & \hspace{10pt}
2w_3(t)=\alpha_3e^{ia_3t}.
\end{array}$
\vspace{8pt}

\noindent Hence, if we let $\beta_1,\beta_1,\beta_3:\R\rightarrow\C$
be differentiable functions such that

\vspace{8pt}
$\begin{array}{rrr}
\hspace{6pt}p_1(t)=ie^{ia_1t}\beta_1(t), & \hspace{15pt} p_2(t)=e^{ia_2t}\beta_2(t), & \hspace{14pt} p_3(t)=e^{ia_3t}\beta_3(t),
\end{array}$
\vspace{8pt}

\noindent we have the following result.

\begin{prop}
\label{matrix1} Using the notation above, \eq{9.2.1}-\eq{9.2.4} can
be written as the following matrix equation for the functions
$x,\beta_1,\beta_2,\beta_3$:
\begin{equation*}
\begin{split}
\frac{d}{dt}\left(\begin{array}{c}
x \\
\beta_1 \\
\beta_2 \\
\beta_3 \\
\bar{\beta}_1 \\
\bar{\beta}_2 \\
\bar{\beta}_3
\end{array} \right)\! = \frac{i}{2}\left(\begin{array}{rrrrrrr}
0 & -\frac{\alpha_1}{2} & \frac{\alpha_2}{2} &
\frac{\alpha_3}{2} & \frac{\alpha_1}{2} &
-\frac{\alpha_2}{2} & -\frac{\alpha_3}{2} \\
\alpha_1 & -2a_1 & 0 & 0 & 0 & -\alpha_3 & -\alpha_2 \\
\alpha_2 & 0 & -2a_2 & 0 & \alpha_3 & 0 & \alpha_1 \\
\alpha_3 & 0 & 0 & -2a_3 & \alpha_2 & \alpha_1 & 0 \\
-\alpha_1 & 0 & \alpha_3 & \alpha_2 & 2a_1 & 0 & 0 \\
-\alpha_2 & -\alpha_3 & 0 & -\alpha_1 & 0 & 2a_2 & 0 \\
-\alpha_3 & -\alpha_2 & -\alpha_1 & 0 & 0 & 0 & 2a_3
\end{array}\right)
\left(\begin{array}{c}
x \\
\beta_1 \\
\beta_2 \\
\beta_3 \\
\bar{\beta}_1 \\
\bar{\beta}_2 \\
\bar{\beta}_3
\end{array} \right)\!\!.
\end{split}
\end{equation*}
\end{prop}

\begin{proof}
Using \eq{9.2.1},
\begin{align*}
\frac{dx}{dt} & = \frac{1}{2}\,{\rm Im}(\alpha_1\beta_1 -\alpha_2
\beta_2 -\alpha_3\beta_3) \\
& =  \frac{i}{4}\left(\alpha_1(\bar{\beta}_1-\beta_1)+\alpha_2
(\beta_2-\bar{\beta}_2)+\alpha_3(\beta_3-\bar{\beta}_3)\right),
\end{align*}
\noindent which gives the first row in the matrix equation above.
Since $a_1+a_2+a_3=0$, equation \eq{9.2.2} for $p_1$ shows that
\begin{equation*}
\begin{split}
i\frac{d\beta_1}{dt}-a_1\beta_1 = \frac{1}{2}(-\alpha_1x+\alpha_2
\bar{\beta}_3+\alpha_3\bar{\beta}_2),
\end{split}
\end{equation*}
\noindent which, upon rearrangement, gives the second row in the
matrix equation above. The calculation of the rest of the rows
follows in a similar fashion.
\end{proof}

\noindent In order to solve the matrix equation given in Proposition
\ref{matrix1}, we find the eigenvalues and corresponding
eigenvectors of the matrix.

\begin{prop}
\label{matrix2} Let $T$ denote the $7\times 7$ real matrix given
in Proposition \ref{matrix1} and let ${\bf
a}=(0,\alpha_1,\alpha_2,\alpha_3,\alpha_1,\alpha_2,\alpha_3)^{\rm
T}$, where ${}^{\rm T}$ denotes transpose.
 Then there exist nonzero vectors ${\bf b}_{\pm}$, ${\bf c}_{\pm}$,
${\bf d}_{\pm}\in\R^7$ such that

\vspace{8pt} $\begin{array}{cccc} T {\bf a}=0, & \hspace{9pt}
T{\bf b}_{\pm} = \pm\lambda{\bf b}_{\pm}, &\hspace{9pt}
 T{\bf c}_{\pm} = \pm\lambda{\bf c}_{\pm}, &  \hspace{9pt}
T{\bf d}_{\pm} = \pm 3\lambda{\bf d}_{\pm},
\end{array}$
\vspace{8pt}

\noindent where $\lambda>0$ is such that $\lambda^2 = a_2^2 -
a_1a_3$ and
\begin{eqnarray}
 {\bf b}_+ = (b_1, b_2, 0, b_3, b_4, 0, b_5)^{{\rm T}}, & &
{\bf b}_- = (b_1, b_4, 0, b_5, b_2, 0, b_3)^{{\rm T}}, \nonumber \\
\label{evectors}
{\bf c}_+ = (c_1, 0, c_2, c_3, 0, c_4, c_5)^{{\rm T}}, &
&{\bf c}_- = (c_1, 0, c_4, c_5, 0, c_2, c_3)^{{\rm T}}, \\
{\bf d}_+ = (0, d_1, d_2, d_3, d_4, d_5, d_6)^{{\rm T}}, &
&{\bf d}_- = (0, d_4, d_5, d_6, d_1, d_2, d_3)^{{\rm T}}, \nonumber
\end{eqnarray}
\noindent for constants
$b_1,\ldots,b_5,c_1,\ldots,c_5,d_1,\ldots,d_6\in\R$. In particular,
the pairs $\{{\bf b}_{\pm},{\bf c}_{\pm}\}$ are linearly
independent.
\end{prop}
\begin{proof}
Most of the results in this proposition are found by direct
calculation using Maple.  The only point to note is that if ${\bf
w}$ is a $\mu$-eigenvector of $T$, for some $\mu\in\R$, and we write
${\bf w} = (\begin{array}{ccc} x & {\bf y} & {\bf z}
\end{array})^{{\rm T}}$, where $x\in\R$ and ${\bf y},{\bf
z}\in\R^3$, then $ \tilde{{\bf w}} = (\begin{array}{ccc} x &  {\bf
z} & {\bf y} \end{array})^{{\rm T}}$ is a $-\mu$-eigenvector of $T$,
and hence we can cast the eigenvectors of $T$ into the form as given
in \eq{evectors}.
\end{proof}

\noindent From this result we can write down the general solution to
the matrix equation given in Proposition \ref{matrix1}:
\begin{align}
&\left(\begin{array}{c}\label{gensoln} x \\
\beta_1 \\
\beta_2 \\
\beta_3 \\
\bar{\beta}_1 \\
\bar{\beta}_2 \\
\bar{\beta}_3 \end{array}\right)  =
A\left(\begin{array}{c} 0 \\ \alpha_1 \\ \alpha_2 \\ \alpha_3 \\ \alpha_1 \\
\alpha_2 \\ \alpha_3 \end{array}\right) +
B_+ e^{\frac{i}{2}\lambda t}\left(\begin{array}{c}
b_1 \\ b_2 \\ 0 \\ b_3 \\ b_4 \\ 0 \\ b_5 \end{array}\right) +
B_- e^{-\frac{i}{2}\lambda t}\left(\begin{array}{c}
b_1 \\ b_4 \\ 0 \\ b_5 \\ b_2 \\ 0 \\ b_3 \end{array}\right)
\\
 &{}\!+ C_+ e^{\frac{i}{2}\lambda t}\!\left(\begin{array}{c} c_1 \\
0
\\ c_2 \\ c_3 \\ 0 \\ c_4 \\ c_5 \end{array}\right)\! + C_-
e^{-\frac{i}{2}\lambda t}\!\left(\begin{array}{c} c_1 \\ 0 \\ c_4 \\
c_5 \\ 0 \\ c_2 \\ c_3 \end{array}\right)\! +\! D_+
e^{\frac{3i}{2}\lambda t}\!\left(\begin{array}{c} 0 \\ d_1 \\ d_2 \\
d_3 \\ d_4 \\ d_5 \\ d_6 \end{array}\right)\! +\! D_-
e^{-\frac{3i}{2}\lambda t}\!\left(\begin{array}{c} 0 \\ d_4 \\ d_5
\\ d_6 \\ d_1 \\ d_2 \\ d_3 \end{array}\right)\nonumber
\end{align}
for constants $A,B_{\pm},C_{\pm},D_{\pm}\in\C$.  However, the last
three rows in this equation are equal to the complex conjugate of
the three rows above them, which implies that $B_-=\bar{B}_+$,
$C_-=\bar{C}_+$, and $D_-=\bar{D}_+$.  Moreover, if we translate
$\R^2$, as given in the evolution data, from $(y_1,y_2)$ to $(y_1-A,
y_2)$, then ${\bf w}_j$ is unaltered for $j=1,2,3$ but ${\bf w}_4$
is mapped to ${\bf w}_4 - A{\bf w}_1$. Therefore, we can set $A=0$.

\medskip

From the discussion above, we may now write down the general
solution to \eq{9.2.1}-\eq{9.2.4} and \eq{9.2.5}-\eq{9.2.8} and then
simply integrate equations \eq{9.2.9}-\eq{9.2.12} to give an
explicit description of some associative 3-folds constructed using
our second evolution equation.  This result is given below.

\begin{thm}
\label{mainresult}
Define functions $x,y,z:\R\rightarrow\R$ and $w_j,p_j,q_j,r_j:\R\rightarrow\C$
for $j=1,2,3$ by:
$$ 2w_1(t) = i\alpha_1e^{ia_1t}, \qquad
2w_2(t)=\alpha_2e^{ia_2t} ,  \qquad 2w_3(t)=\alpha_3e^{ia_3t},
$$

\noindent where $\alpha_1,\alpha_2,\alpha_3$ are positive constants
such that $\alpha_1^{-2}=\alpha_2^{-2}+\alpha_3^{-2}$ and
$a_1,a_2,a_3$ are given in \eq{as};
\begin{align*}
x(t) & =  2\hspace{2pt} {\rm Re}\left(Bb_1e^{\frac{i}{2}\lambda t}+Cc_1e^{\frac{i}{2}\lambda t}\right), \\
p_1(t) & =  ie^{ia_1t}\left(Bb_2e^{\frac{i}{2}\lambda t}+\bar{B}b_4e^{-\frac{i}{2}\lambda t}+Dd_1e^{\frac{3i}{2}\lambda t}+\bar{D}d_4e^{-\frac{3i}{2}\lambda t}\right), \\
p_2(t) & =  e^{ia_2t}\left(Cc_2e^{\frac{i}{2}\lambda t}+\bar{C}c_4e^{-\frac{i}{2}\lambda t}+Dd_2e^{\frac{3i}{2}\lambda t}+\bar{D}d_5e^{-\frac{3i}{2}\lambda t}\right), \\
p_3(t) & =  e^{ia_3t}\Big((Bb_3+Cc_3)e^{\frac{i}{2}\lambda t}\!+\!(\bar{B}b_5+\bar{C}c_5)e^{-\frac{i}{2}\lambda t}\!+\!Dd_3e^{\frac{3i}{2}\lambda t}\!+\!\bar{D}d_6e^{-\frac{3i}{2}\lambda t}\Big), \\
y(t) & =  2\hspace{2pt} {\rm Re}\left(B^{\prime}b_1e^{\frac{i}{2}\lambda t}+C^{\prime}c_1e^{\frac{i}{2}\lambda t}\right), \\
q_1(t) & =  ie^{ia_1t}\left(B^{\prime}b_2e^{\frac{i}{2}\lambda t}\!+\bar{B^{\prime}}b_4e^{-\frac{i}{2}\lambda t}\!+D^{\prime}d_1e^{\frac{3i}{2}\lambda t}\!+\bar{D^{\prime}}d_4e^{-\frac{3i}{2}\lambda t}\right), \\
q_2(t) & =  e^{ia_2t}\left(C^{\prime}c_2e^{\frac{i}{2}\lambda t}+\bar{C^{\prime}}c_4e^{-\frac{i}{2}\lambda t}+D^{\prime}d_2e^{\frac{3i}{2}\lambda t}+\bar{D^{\prime}}d_5e^{-\frac{3i}{2}\lambda t} \right), \\
q_3(t) & =  e^{ia_3t}\Big(\!(B^{\prime}b_3\!+\!C^{\prime}c_3)e^{\frac{i}{2}\lambda t}\!\!+\!(\bar{B^{\prime}}b_5\!+\!\bar{C^{\prime}}c_5)e^{-\frac{i}{2}\lambda t}\!\!+\!D^{\prime}d_3e^{\frac{3i}{2}\lambda t}\!+\!\bar{D^{\prime}}d_6e^{-\frac{3i}{2}\lambda t}\Big), \\
\frac{dz}{dt} & =  {\rm Im}(\bar{p}_1q_1 - \bar{p}_2q_2 - \bar{p}_3q_3), \\
\frac{dr_1}{dt} & =  ixp_1 + iyq_1 +\overline{p_2p_3} + \overline{q_2q_3}, \\
\frac{dr_2}{dt} & =  ixp_2 - iyq_2 -\overline{p_3p_1} + \overline{q_3q_1}, \\
\frac{dr_3}{dt} & =  ixq_3 + iyp_3 -\overline{p_1q_2} -
\overline{p_2q_1},
\end{align*}
\noindent where the real constants $\lambda$ and $b_j,c_j,d_j$ are
as defined in Proposition \ref{matrix2} and $B,B^{\prime},C,
C^{\prime},D,D^{\prime}\in\C$ are arbitrary constants.

Define a subset $M$ of $\R^7=\R\oplus\C^3$ by:
\begin{align}
M  = \bigg\{ \bigg(&y_1y(t)-y_2x(t)+z(t), \hspace{2pt} \frac{1}{2}(y_1^2+y_2^2)w_1(t)+y_1p_1(t)+y_2q_1(t)+r_1(t), \nonumber\\
&\frac{1}{2}(y_1^2-y_2^2)w_2(t)+y_1p_2(t)-y_2q_2(t)+r_2(t), \nonumber\\
\label{affM}&y_1y_2w_3(t)+y_1q_3(t)+y_2p_3(t)+r_3(t)\bigg):y_1,y_2,t\in\R\bigg\}.
\end{align}
Then $M$ is an associative 3-fold in $\R^7$.
\end{thm}

\noindent We now count parameters for the associative 3-folds
constructed by
 Theorem \ref{mainresult}.  There are four real parameters ($\alpha_1,\alpha_2,\alpha_3$
  and the constant of integration for $z(t)$)
and nine complex parameters ($B,B^{\prime},C,C^{\prime},D,
D^{\prime}$ and the three constants of integration for
$r_1(t),r_2(t), r_3(t)$), which makes a total of 22 real parameters.
The relationship $\alpha_1^{-2}=\alpha_2^{-2}+\alpha_3^{-2}$ then
reduces the number of parameters by one to 21.  Recall that we have
the symmetry groups $\GL(2,\R)\ltimes\R^2$ and
$\text{G}_2\ltimes\R^7$ for these associative 3-folds.  By the
arguments proceeding Theorem \ref{evolvethm} and the proof of
\cite[Proposition 9.1]{Joy4}, we have used the freedom in
$\text{G}_2$ transformations and rotations in $\GL(2,\R)$ to
transform our solutions $\textbf{w}_1,\textbf{w}_2,\textbf{w}_3$ of
\eq{ch5.1}-\eq{ch5.3} to solutions of \eq{Joyeq1}-\eq{Joyeq3},
satisfying \eq{Joycond1}, of the form $\textbf{w}_1=(0,w_1,0,0)$,
$\textbf{w}_2=(0,0,w_2,0)$, $\textbf{w}_3=(0,0,0,w_3)$.  We have
also used translations in $\R^2$ to set the constant $A$ in
\eq{gensoln} and the corresponding constant $A^{\prime}$ in the
general solution to \eq{9.2.5}-\eq{9.2.8} both to zero. Therefore,
the remaining symmetries available are dilations in $\GL(2,\R)$ and
translations in $\R^7$, which reduce the number of parameters by
eight to 13.  Translation in time, say $t\mapsto t+t_0$, corresponds
to multiplying $B,B^{\prime},C,C^{\prime}$ by $e^{\frac{i}{2}\lambda
t_0}$ and $D, D^{\prime}$ by $e^{\frac{3i}{2}\lambda t_0}$, which
thus lowers the parameter count by one. We conclude that the
dimension of the family of associative 3-folds generated by Theorem
\ref{mainresult} is 12, whereas the dimension of the whole family
generated by Theorem \ref{evolvethm} is 14.

\subsection{Periodicity}

Note that the solutions to Theorem \ref{mainresult} are all linear
combinations of terms of the form $e^{i(a_j+m\lambda)t}$ for
$j=1,2,3$ and $m=0,\pm \frac{1}{2},\pm 1,\pm \frac{3}{2},\pm 2,\pm
3$, since $a_j\pm n \lambda\neq 0$ for $n=0,1,2,3$, which ensures
that $r_1,r_2,r_3$ do not have any linear terms in $t$.  It is
therefore reasonable to search for associative 3-folds $M$ as in
\eq{affM} that are \emph{periodic} in $t$.  Define a map
$F:\R^3\rightarrow\R^7$ by \eq{singmap}, so that
$M=\text{Image}\,F$. Then $M$ is periodic if and only if there
exists some constant $T>0$ such that $F(y_1,y_2,t+T)=F(y_1,y_2,t)$
for all $y_1,y_2,t\in\R$.

From above, the periods of the exponentials in the functions defined
in Theorem \ref{mainresult} are proportional to
$(a_j+m\lambda)^{-1}$ for $j=1,2,3$ and the values of $m$ given
above. In general $F$ will be periodic if and only if these periods
have a common multiple.  By the definition of the constants $a_j$,
we can write $a_2=-xa_1$ and $a_3=(x-1)a_1$ for some $x\in (0,1)$.
Then $\lambda^2= a_2^2-a_1a_3=a_1^2(x^2-x+1)$ and, if we let
$y=\sqrt{x^2-x+1}$, we deduce that $\lambda=-ya_1$ since $a_1<0$ and
$\lambda,y>0$. The periods thus have a common multiple if and only
if $x$ and $y$ are rational.  We have therefore reduced the problem
to finding rational points on the conic $y^2=x^2-x+1$. This is a
standard problem in number theory and is identical to the one solved
by Joyce \cite[$\S$11.2]{Joy4}, so we are able to prove the
following result.

\begin{thm}
\label{periodicthm} Given $s\in (0,\frac{1}{2})\cap\Q$, Theorem
\ref{mainresult} gives a family of closed associative 3-folds in
$\R^7$ whose generic members are nonsingular immersed\ 3-folds
diffeomorphic to $\mathcal{S}^1\times\R^2$.
\end{thm}

\begin{proof}
Let $s\in (0,\frac{1}{2})\cap\Q$ and write $s=\frac{p} {q}$ where
$p,q$ are coprime positive integers.  Then, by the work in
\cite[p.390]{Joy4}, we define $a_1,a_2,a_3,\lambda$ either by
\begin{align*}
a_1&=p^2-q^2,&a_2&=q^2-2pq,& a_3&=2pq-p^2, & \lambda&=p^2-pq+q^2;\\
\intertext{or, if $p+q$ is divisible by 3, by}
3a_1&=p^2-q^2,&3a_2&=q^2-2pq,&3a_3&=2pq-p^2,&3\lambda&=p^2-pq+q^2.
\end{align*}




\noindent In both cases,
$\text{hcf}\,(a_1,a_2,a_3)=\text{hcf}\,(a_1,a_2,a_3, \lambda)=1$. We
also note that $\lambda$ is odd since at least one of $p,q$ is odd.
Thus $a_j+m\lambda$ is an integer for integer values of $m$ and half
an integer, but not an integer, for non-integer values of $m$.
Hence, by the form of the functions given in Theorem
\ref{mainresult} and equation \eq{singmap} for $F$,
$F(y_1,y_2,t+2\pi)=F(-y_1,-y_2,t)$ for all $y_1,y_2,t$.  We deduce
that $F$ has period $4\pi$, using the condition that
hcf$(a_1,a_2,a_3)=1$.

If we define an action of $\Z$ on $\R^3$ by requiring,  for
$n\in\Z$, that $(y_1,y_2,t)$ maps to $\left((-1)^ny_1,
(-1)^ny_2,t+2n\pi\right)$, then we can consider $F$ as a map from
the quotient of $\R^3$ by $\Z$ under this action. Since this
quotient is diffeomorphic to $\mathcal{S}^1\times\R^2$ and
generically $F$ is an immersion, $M=\text{Image}\,F$ is generically
an immersed 3-fold diffeomorphic to $\mathcal{S}^1\times \R^2$.
\end{proof}

 Joyce \cite{Joy4} has considered the asymptotic behaviour
of the SL 3-folds constructed by Theorem \ref{Joythm}(iv) at
infinity, which is dependent on the quadratic terms in $F$.
However, since solutions ${\bf w}_1,{\bf w}_2,{\bf w}_3$ in
Theorem \ref{evolvethm} are essentially equivalent to solutions
${\bf z}_1,{\bf z}_2,{\bf z}_3$ in Theorem \ref{Joythm}, the
asymptotic behaviour of the 3-folds given by Theorem
\ref{periodicthm} must be identical to that found by Joyce
\cite[p.391]{Joy4}.  We first make a definition and then state our
result.

\begin{dfn}\label{wasym}
Let $M,M_0$ be closed $m$-dimensional submanifolds of $\R^n$ and let
$k<1$. We say that $M$ is \emph{asymptotic with order $O(r^k)$ at
infinity in $\R^n$ to $M_0$} if there exist $R>0$, some compact
subset $K$ of $M$ and a diffeomorphism
$\Phi:M_0\setminus\bar{B}_R\rightarrow M\setminus K$ such that
$$|\Phi(\textbf{x})-\textbf{x}|=O(r^k)\quad\text{as $r\rightarrow\infty$},$$
where $r$ is the radius function on $\R^n$ and $\bar{B}_R$ is the
closed ball of radius $R$.
\end{dfn}

\begin{thm}
\label{diverge} Every closed associative 3-fold defined by $s\in
(0,\frac{1}{2})\cap\Q$, as given in Theorem \ref{periodicthm}, is
asymptotic with order $O(r^{\frac{1}{2}})$ at infinity in $\R^7$
to a double cover of the SL $T^2$ cone defined by:
\begin{equation*}
\left\{(0,ie^{ia_1t}x_1, e^{ia_2t}x_2, e^{ia_3t}x_3):
x_1,x_2,x_3,t\in\R,\, x_1\geq 0,\, \sum_{i=1}^3a_ix_i^2=0\right\}
\end{equation*}
\noindent where the constants $a_1,a_2,a_3$ are defined by $s$ as in the
proof of Theorem \ref{periodicthm}.
\end{thm}

\noindent The associative 3-folds in Theorem \ref{periodicthm}
actually diverge away from the SL cone given above, but Theorem
\ref{diverge} gives a measure of the rate of divergence.

We now show that if an associative 3-fold  $M$ were to converge to
an SL 3-fold at infinity, then $M$ would in fact be SL, which we
know is not the case for generic members of the family given by
Theorem \ref{periodicthm}.

\begin{thm}
\label{SLconverge} Suppose $M$ is an associative 3-fold in
$\R^7=\R\oplus\C^3$ and that $L$ is an SL 3-fold in $\C^3$.
Suppose further that $M$ is asymptotic with order $O(r^k)$ at
infinity in $\R^7$ to $L$, where $k<0$.  Then $M$ is an SL 3-fold
in $\C^3$ embedded in $\R^7$.
\end{thm}

\noindent To prove Theorem \ref{SLconverge} we need two results. The
first is a \emph{maximum principle} for harmonic functions due to
Hopf \cite[p. 12]{Lawson}.

\begin{thm}
\label{harmonic1} Let $f$ be a smooth function on a Riemannian
manifold $M$. Suppose $f$ is harmonic, i.e. $d^{\ast}df=0$ where
$d^{\ast}$ is the formal adjoint of $d$.  If $f$ assumes a local
maximum (or minimum) at a point in $M\setminus\partial M$ then $f$
is constant.
\end{thm}

\noindent The second is an elementary result from the theory of minimal
submanifolds \cite[Corollary 9]{Lawson}.

\begin{thm}
\label{harmonic2} Let $M$ be a submanifold of\/ $\R^n$, for some
$n$, with immersion $\iota$. Then $M$ is a minimal submanifold if
and only if $\iota$ is harmonic.
\end{thm}

\noindent Here the function $\iota:M\rightarrow\R^n$ is harmonic if
and only if each of the components of $\iota$ mapping to $\R$ is
harmonic.  We now prove Theorem \ref{SLconverge}.

\begin{proof}[Theorem \ref{SLconverge}]
Since $M$ is an associative 3-fold in $\R^7$, $M$ is a minimal
submanifold of $\R^7$ \cite[Theorem II.4.2]{HarLaw}. Therefore, the
embedding of $M$ in $\R^7$ is harmonic by Theorem \ref{harmonic2}.
In particular, if we write coordinates on $M$ as $(x_1,\ldots,x_7)$,
$x_1$ is harmonic.  We may assume, without loss of generality, that
the SL 3-fold $L$ to which $M$ converges lies in
$\{0\}\times\C^3\subseteq\R^7$. Since $M$ is asymptotic to $L$ at
infinity with order $O(r^k)$ where $k<0$, $x_1\rightarrow 0$ as
$r\rightarrow\infty$.  Suppose $x_1$ is not identically zero. Then
$x_1$ assumes a strict maximum or minimum at some point in the
interior of $M$.  By Theorem \ref{harmonic1}, $x_1$ is therefore
constant, which contradicts the assumption that $x_1$ was not
identically zero. Hence $x_1\equiv 0$ and $M$ is an SL 3-fold in
$\C^3$.
\end{proof}

\section{Ruled Associative 3-folds}
\label{rul}
\setcounter{thm}{0}

In this final section we focus on \emph{ruled} 3-folds and apply our
ideas of evolution equations to give methods for constructing
associative examples. This is a generalisation of the work in
Joyce's paper \cite{Joy5} on ruled SL 3-folds in $\C^3$ and it is
from this source that we take the definitions below.  By a
\emph{cone} in $\R^7$ we shall mean a submanifold of $\R^7$ which is
invariant under dilations and is nonsingular except possibly at 0. A
cone $C$ is said to be \emph{two-sided} if $C=-C$.

\begin{dfn}
\label{ruled} Let $M$ be a 3-dimensional submanifold of $\R^7$.  A
\emph{ruling} of $M$ is a pair $(\Sigma,\pi)$, where $\Sigma$ is a
2-dimensional manifold and $\pi:M \rightarrow\Sigma$ is a smooth
map, such that for all $\sigma\in\Sigma$ there exist ${\bf
v}_{\sigma}\in \mathcal{S}^6$, ${\bf w}_{\sigma}\in\R^7$ such that
$\pi^{-1}(\sigma) =\{r{\bf v}_{\sigma}+{\bf w}_{\sigma}:r\in\R\}$.
Then the triple $(M,\Sigma, \pi)$ is a \emph{ruled submanifold} of
$\R^7$.

\noindent An \emph{r-orientation} for a ruling $(\Sigma,\pi)$ of $M$
is a choice of orientation for the affine straight line
$\pi^{-1}(\sigma)$ in $\R^7$, for each $\sigma\in\Sigma$, which
varies smoothly with $\sigma$.  A ruled submanifold with an
r-orientation for the ruling is called an \emph{r-oriented ruled
submanifold}.

\noindent Let $(M,\Sigma,\pi)$ be an r-oriented ruled submanifold.
For each $\sigma \in\Sigma$, let $\phi(\sigma)$ be the unique unit
vector in $\R^7$ parallel to $\pi^{-1}(\sigma)$ and in the positive
direction with respect to the orientation on $\pi^{-1}(\sigma)$,
given by the r-orientation. Then $\phi:\Sigma\rightarrow
\mathcal{S}^6$ is a smooth map. Define $\psi:\Sigma \rightarrow\R^7$
such that, for all $\sigma\in\Sigma$, $\psi(\sigma)$ is the unique
vector in $\pi^{-1}(\sigma)$ orthogonal to $\phi(\sigma)$.  Then
$\psi$ is a smooth map and we may write:
\begin{equation}
\begin{split}
\label{ruled1}
M=\{r\phi(\sigma)+\psi(\sigma):\sigma\in\Sigma,\hspace{2pt} r\in\R\}.
\end{split}
\end{equation}
Define the \emph{asymptotic cone} $M_0$ of a ruled submanifold $M$
by:
\begin{equation*}
\begin{split}
M_0=\{{\bf v}\in\R^7:\text{${\bf v}$ is parallel to
$\pi^{-1}(\sigma)$ for some $\sigma\in\Sigma$}\}.
\end{split}
\end{equation*}
If $M$ is also r-oriented then
\begin{equation}
\begin{split}
\label{ruled2}
M_0=\{r\phi(\sigma):\sigma\in\Sigma,\hspace{2pt}r\in\R\}
\end{split}
\end{equation}
\noindent and is usually a 3-dimensional two-sided cone; that is,
whenever $\phi$ is an immersion.
\end{dfn}

Note that we can consider any r-oriented ruled submanifold as being
defined by two maps $\phi,\psi$ as given in Definition \ref{ruled}.
Hence, r-oriented ruled associative 3-folds may be constructed by
evolution equations for $\phi,\psi$.

\medskip

Suppose we have a 3-dimensional two-sided cone $M_0$ in $\R^7$. The
\emph{link} of $M_0$, $M_0\cap \mathcal{S}^6$, is a nonsingular
2-dimensional submanifold of $\mathcal{S}^6$ closed under the action
of $-1:\mathcal{S}^6 \rightarrow \mathcal{S}^6$. Let $\Sigma$ be the
quotient of the link by the $\pm1$ maps on $\mathcal{S}^6$. Clearly,
$\Sigma$ is a nonsingular 2-dimensional manifold. Define
$\tilde{M_0} \subseteq\Sigma\times\R^7$ by:
\begin{equation*}
\begin{split}
\tilde{M_0}=\{(\{\pm \sigma\},r\sigma):\sigma\in M_0\cap
\mathcal{S}^6,\hspace{2pt}r\in\R\}.
\end{split}
\end{equation*}
\noindent Then $\tilde{M_0}$ is a nonsingular 3-fold.  Define
$\pi:\tilde{M_0}\rightarrow\Sigma$ by $\pi(\{\pm
\sigma\},r\sigma)=\{\pm\sigma\}$ and
$\iota:\tilde{M_0}\rightarrow\R^7$ by $\iota(\{\pm
\sigma\},r\sigma)= r\sigma$.  Note that $\iota(\tilde{M_0})=M_0$ and
that $\iota$ is an immersion except on $\iota^{-1}(0)\cong\Sigma$,
so we may consider $\tilde{M_0}$
 as a singular immersed submanifold of $\R^7$. Hence $(\tilde{M_0},
\Sigma,\pi)$ is a ruled submanifold of $\R^7$.  Therefore, we can
regard $M_0$ as a ruled submanifold and dispense with $\tilde{M_0}$.
Suppose further that $M_0$ is an r-oriented two-sided cone. We can
thus write $M_0$ in the form \eq{ruled1} for maps $\phi,\psi$, as
given in Definition
 \ref{ruled}, and see that $\psi$ must be identically zero.  It is also clear that any
ruled submanifold defined by $\phi,\psi$ with $\psi\equiv0$ is an
r-oriented two-sided cone.

We now justify the terminology of asymptotic cone as given in
Definition \ref{ruled}. For this, we need to define the term
\emph{asymptotically conical with order $O(r^{\alpha})$}, where $r$
is the radius function on $\R^7$.

\begin{dfn}
\label{asym} Let $M_0$ be a closed cone in $\R^7$ and let $M$ be a
closed nonsingular submanifold in $\R^7$.  We say that $M$ is
\emph{asymptotically conical to $M_0$ with order $O(r^{\alpha})$},
for some $\alpha<1$, if there exist some constant $R>0$, a compact
subset $K$ of $M$ and a diffeomorphism
$\Phi:M_0\setminus\bar{B}_R\rightarrow M\setminus K$ such that
\begin{equation}
\label{asymconds} |\nabla^{k}(\Phi({\bf x})-I({\bf
x}))|=O(r^{\alpha-k}) \quad \text{for $k=0,1,2,\ldots$ as
$r\rightarrow\infty$}
\end{equation}
\noindent where $\bar{B}_R$ is the closed ball of radius $R$ in
$\R^7$ and $I:M_0\rightarrow\R^7$ is the inclusion map.  Here
$|\,.\,|$ is calculated using the cone metric on
$M_0\setminus\bar{B}_R$, and $\nabla$ is a combination of the
Levi--Civita connection derived from the cone metric and the flat
connection on $\R^n$, which acts as partial differentiation.
\end{dfn}

\noindent Suppose that $M$ is an r-oriented ruled submanifold and
let $M_0$ be its asymptotic cone.  Writing $M$ in the form
\eq{ruled1} and $M_0$ in the form \eq{ruled2} for maps $\phi,\psi$,
define a diffeomorphism $\Phi:M_0\setminus \bar{B}_1\rightarrow
M\setminus K$, where $K$ is some compact subset of $M$, by
$\Phi(r\phi(\sigma))=r\phi(\sigma)+\psi(\sigma)$ for all
$\sigma\in\Sigma$ and $|r|>1$.  If $\Sigma$ is compact, so that
$\psi$ is bounded, then $\Phi$ satisfies \eq{asymconds} as given in
Definition \ref{asym} for $\alpha=0$, which shows that $M$ is
asymptotically conical to $M_0$ with order $O(1)$.




\subsection{The associative condition}


Let $\Sigma$ be a 2-dimensional, connected, real analytic manifold,
let $\phi:\Sigma\rightarrow \mathcal{S}^6$ be a real analytic
immersion and let $\psi:\Sigma \rightarrow\R^7$ be a real analytic
map.  Define $M$ by \eq{ruled1}, so that $M$ is the image of the map
$\iota:\R\times\Sigma\rightarrow\R^7$ given by
$\iota(r,\sigma)=r\phi(\sigma)+\psi(\sigma)$. Clearly,
$\R\times\Sigma$ is an r-oriented ruled submanifold with ruling
$(\Sigma,\pi)$, where $\pi$ is given by $\pi(r,\sigma)=\sigma$.
Since $\phi$ is an immersion, $\iota$ is an immersion almost
everywhere in $\R\times\Sigma$ and thus $M$ is an r-oriented ruled
submanifold.

We now suppose that $M$ is associative in order to discover the
conditions that this imposes upon $\phi,\psi$.  Note that the
asymptotic cone $M_0$ of $M$, given by $\eq{ruled2}$, is the image
of $\R\times\Sigma$ under the map $\iota_0$, defined by
$\iota_0(r,\sigma)=r\phi( \sigma)$. Since $\phi$ is an immersion,
$\iota_0$ is an immersion except at $r=0$, so $M_0$ is a
3-dimensional cone which is nonsingular except at 0.


Let $p\in M$. There exist $r\in\R$, $\sigma\in\Sigma$ such that
$p=r\phi(\sigma)+\psi(\sigma)$.  Choose local coordinates $(s,t)$
near $\sigma$ in $\Sigma$. Then $T_pM=\langle x,y,z \rangle_{\R}$,
where $x=\phi(\sigma)$, $y=r\frac{\partial \phi}{\partial
s}(\sigma)+ \frac{\partial \psi}{\partial s}(\sigma)$ and
$z=r\frac{\partial \phi}{\partial
 t}(\sigma)+\frac{\partial \psi}{\partial t}(\sigma)$.  Since $M$ is
 associative, $T_pM$ is an associative 3-plane, which by Proposition \ref{ass}
occurs if and only if $[x,y,z]=0$. This condition forces a quadratic
in $r$ to vanish, and thus the coefficient of each power of $r$ must
be zero as this condition should hold for all $r\in\R$.  The
following set of equations must therefore hold in $\Sigma$:
\begin{align}
\label{assruled1}
\left[\phi,\frac{\partial\phi}{\partial s},
\frac{\partial\phi}{\partial t}\right] &= 0, \\
\label{assruled2}
\left[\phi,\frac{\partial\phi}{\partial s},
\frac{\partial\psi}{\partial t}\right] +
\left[\phi,\frac{\partial\psi}{\partial s},
\frac{\partial\phi}{\partial t}\right] & =  0, \\
\label{assruled3} \left[\phi,\frac{\partial\psi}{\partial s},
\frac{\partial\psi}{\partial t}\right] &= 0.
\end{align}
Note firstly that, if we do not suppose $M$ to be associative but
that \eq{assruled1}-\eq{assruled3} hold locally in $\Sigma$, then
following the argument above we see that each tangent space to $M$
must be associative and hence that $M$ is associative.  Moreover,
\eq{assruled1} is equivalent to having that tangent spaces to points
 of the form $r\phi(\sigma)$, for $r\in\R,\sigma\in\Sigma$, are associative,
which is precisely the condition for the asymptotic cone $M_0$ to be
associative.  We may therefore deduce the following result.

\begin{prop}
\label{asscone} Let $M$ be an r-oriented ruled associative 3-fold in
$\R^7$ and let $M_0$ be the asymptotic
 cone of $M$.  Then $M_0$ is an associative cone in $\R^7$ provided
 it is 3-dimensional.
\end{prop}

Since $M_0$ is associative, $\varphi$ is a non-vanishing 3-form on
$M_0$ that defines the orientation on $M_0$.  This forces $\Sigma$
to be oriented, for if $(s,t)$ are some local coordinates on
$\Sigma$, then we can define them to be oriented by imposing the
condition that
\begin{equation*}
\begin{split}
\varphi\left(\phi,\frac{\partial\phi}{\partial
s},\frac{\partial\phi} {\partial t}\right)
>0.
\end{split}
\end{equation*}
\noindent In addition, if $g$ is the natural metric on
$\mathcal{S}^6$, then the pullback $\phi^*(g)$ is a metric on
$\Sigma$ making it a \emph{Riemannian} 2-fold, since
$\phi:\Sigma\rightarrow \mathcal{S}^6$ is an immersion. Therefore we
can consider $\Sigma$ as an oriented Riemannian 2-fold and hence it
has a natural \emph{complex structure}, which we denote as $J$.
Locally in $\Sigma$ we can choose a holomorphic coordinate $u=s+it$,
and so the corresponding real coordinates $(s,t)$ satisfy the
condition $J(\frac{\partial}{\partial s})=\frac{\partial}{\partial
t}$. Following Joyce \cite[p.241]{Joy5}, we say that local real
coordinates $(s,t)$ on $\Sigma$ that have this property are
\emph{oriented conformal coordinates}.

We now use oriented conformal coordinates in the proof of the next
result, which gives neater equations for maps $\phi,\psi$ defining
an r-oriented ruled associative 3-fold.

\begin{thm}
\label{ruledthm} Let $\Sigma$ be a connected real analytic 2-fold,
let $\phi:\Sigma\rightarrow \mathcal{S}^6$ be a real analytic
immersion and let $\psi:\Sigma \rightarrow \R^7$ be a real analytic
map.  Let $M$ be defined by \eq{ruled1}.  Then $M$ is associative if
and only if
\begin{equation}
\begin{split}
\label{assruled4} \frac{\partial \phi}{\partial t} =
\phi\times\frac{\partial \phi}{\partial s}
\end{split}
\end{equation}
\noindent and $\psi$ satisfies

\vspace{4pt}
\begin{rlist}
\item
$\frac{\partial\psi}{\partial t}=\phi\times\frac{\partial\psi}{\partial s} +
f\phi$ for some real analytic function $f:\Sigma\rightarrow\R$,

\vspace{4pt}
or
\vspace{4pt}

\item $\frac{\partial\psi}{\partial
s}(\sigma),\frac{\partial\psi}{\partial t}
(\sigma)\in\langle\phi(\sigma),\frac{\partial\phi}{\partial
s}(\sigma), \frac{\partial\phi}{\partial t}(\sigma)\rangle_{\R}$
for all $\sigma\in\Sigma$,
\end{rlist}
\vspace{4pt}

\noindent where $\times$ is defined by \eq{cross2} and $(s,t)$ are oriented
conformal coordinates on $\Sigma$.
\end{thm}

\begin{proof}
Above we noted that \eq{assruled1}-\eq{assruled3} were equivalent to
the condition that $M$ is associative, so we show that
\eq{assruled4} is equivalent to \eq{assruled1} and that (i) and (ii)
are equivalent to \eq{assruled2} and \eq{assruled3}.

Let $\sigma\in\Sigma$, $C=|\frac{\partial\phi}{\partial
s}(\sigma)|>0$.  Since $\phi$ maps to the unit sphere in $\R^7$,
$\phi(\sigma)$ is orthogonal to $\frac{\partial\phi}{\partial
s}(\sigma)$ and $\frac{\partial\phi}{\partial t}(\sigma)$. As
$(s,t)$ are oriented conformal coordinates, we also see that
$\frac{\partial\phi}{\partial s}(\sigma)$ and
$\frac{\partial\phi}{\partial t}(\sigma)$ are orthogonal and that
$|\frac{\partial\phi}{\partial t}(\sigma)|=C$.  We conclude that the
triple $(\phi(\sigma), C^{-1}\frac{\partial\phi}{\partial
s}(\sigma), C^{-1}\frac{\partial\phi}{\partial t}(\sigma))$ is an
oriented orthonormal triad in $\R^7$, and it is the basis for an
associative 3-plane in $\R^7$ if and only if \eq{assruled1} holds at
$\sigma$. Since $\text{G}_2$ acts transitively on the set of
associative 3-planes \cite[Theorem IV.1.8]{HarLaw}, if
\eq{assruled1} holds at $\sigma$ then we can transform coordinates
on $\R^7$ using $\text{G}_2$ so that
\begin{equation*} \phi(\sigma)=e_1, \qquad
\frac{\partial\phi}{\partial s}\,(\sigma)=Ce_2, \qquad
\frac{\partial\phi}{\partial t}\,(\sigma)=Ce_3,
\end{equation*}
where $\{e_1,\ldots,e_7\}$ is a basis for Im $\O\cong\R^7$.  We note
here that \eq{assruled4} holds at $\sigma$ since the cross product
is invariant under $\text{G}_2$. It is clear that, if \eq{assruled4}
holds at $\sigma$, then the 3-plane generated by
$\{\phi(\sigma),\frac{\partial\phi}{\partial
s}(\sigma),\frac{\partial\phi}{\partial t}(\sigma)\}$ is
associative, simply by the definition of the cross product.

Under the change of coordinates of $\R^7$ above, we can write
$\frac{\partial\psi}{\partial s}(\sigma)=a_1e_1+\ldots+a_7e_7$ and
$\frac{\partial\psi}{\partial t}(\sigma)=b_1e_1+\ldots+b_7e_7$ for
real constants $a_j,b_j$ for $j=1,\ldots,7$.  Calculations show that
 $\eq{assruled2}$ holds at $\sigma$ if and only if
\begin{equation}
\label{assruled5}
b_4=-a_5, \hspace{20pt} b_5=a_4, \hspace{20pt} b_6=-a_7, \hspace{20pt} b_7=a_6,
\end{equation}
\noindent and \eq{assruled3} holds at $\sigma$ if and only if
\begin{align}
\label{assruled6.1}
-a_4b_7-a_5b_6+a_6b_5+b_4a_7 &= 0, \\
\label{assruled6.2}
-a_4b_6+a_5b_7+a_6b_4-a_7b_5 &= 0, \\
\label{assruled6.3}
 a_2b_7+a_3b_6-a_6b_3-a_7b_2 &= 0, \\
\label{assruled6.4}
 a_2b_6-a_3b_7-a_6b_2+a_7b_3 &= 0, \\
\label{assruled6.5}
-a_2b_5-a_3b_4+a_4b_3+a_5b_2 &= 0, \\
\label{assruled6.6} -a_2b_4+a_3b_5+a_4b_2-a_5b_3 &= 0.
\end{align}
\noindent Substituting condition \eq{assruled5} into the above
equations, \eq{assruled6.1} and \eq{assruled6.2} are satisfied
immediately and \eq{assruled6.3}-\eq{assruled6.6} become:
\begin{align*}
a_6(a_2-b_3)-a_7(a_3+b_2) &= 0, \\
-a_6(a_3+b_2)-a_7(a_2-b_3) &= 0, \\
-a_4(a_2-b_3)+a_5(a_3+b_2) &= 0, \\
a_4(a_3+b_2)+a_5(a_2-b_3) &= 0.
\end{align*}
\noindent These equations can then be written in matrix form:
\begin{align}
\label{assruled7.1} \left(\begin{array}{rr}
-a_6 & a_7 \\
a_7 & a_6 \end{array}\right)
\left(\begin{array}{c}
a_2-b_3 \\
a_3+b_2 \end{array}\right) & =  0, \\
\label{assruled7.2}
\left(\begin{array}{rr}
-a_4 & a_5 \\
a_5 & a_4 \end{array}\right)
\left(\begin{array}{c}
a_2-b_3 \\
a_3+b_2 \end{array}\right) & =  0.
\end{align}
\noindent We see that equations \eq{assruled7.1} and
\eq{assruled7.2} hold if and only if the vector appearing in both
equations is zero or the determinants of the matrices are zero. We
thus have two conditions which we shall show correspond to (i) and
(ii):
\begin{eqnarray}
\label{assruled8.1}
a_2=b_3, \hspace{20pt} -a_3=b_2; \\
\label{assruled8.2}
a_4=a_5 = 0 =a_6=a_7.
\end{eqnarray}

Using the fact that $\phi(\sigma)=e_1$, \eq{assruled8.1} holds if
and only if
\begin{align*}
\frac{\partial\psi}{\partial t}
(\sigma) &=  b_1e_1-a_3e_2+a_2e_3-a_5e_4+a_4e_5-a_7e_6+a_6e_7 \\
& =  \phi(\sigma)\times\frac{\partial\psi}{\partial
s}(\sigma)+f(\sigma) \phi(\sigma),
\end{align*}
\noindent where $f(\sigma)=b_1$. Therefore, \eq{assruled8.1}
corresponds to condition (i) holding at $\sigma$ by virtue of the
invariance of the cross product under $\text{G}_2$.  The fact that
$f$ is real analytic is immediate from the hypotheses that
$\phi,\psi$ are real analytic and that $\phi$ is nonzero, since
$\phi$ maps to $\mathcal{S}^6$.

Similarly, \eq{assruled8.2} holds if and only if
\begin{equation*}
\frac{\partial\psi}{\partial s}(\sigma)=a_1e_1+a_2e_2+a_3e_3 \hspace{10pt} {\rm and} \hspace{10pt}
\frac{\partial\psi}{\partial t}(\sigma)=b_1e_1+b_2e_2+b_3e_3,
\end{equation*}
\noindent which is equivalent to condition (ii) holding at $\sigma$,
since we may note here that $\langle e_1,e_2,e_3\rangle_{\R}=\langle
\phi(\sigma),\frac{\partial\phi} {\partial
s}(\sigma),\frac{\partial\phi}{\partial t}(\sigma)\rangle_{\R}$.

In conclusion, at each point $\sigma\in\Sigma$, condition (i) or
(ii) holds.  Let $\Sigma_1=\{\sigma\in\Sigma\,:\,\text{(i) holds at
$\sigma$}\}$ and let $\Sigma_2=\{\sigma\in\Sigma\,:\,\text{(ii)
holds at $\sigma$}\}$.  Note that (i) and (ii) are closed conditions
on the real analytic maps $\phi,\psi$.  Therefore, $\Sigma_1$ and
$\Sigma_2$ are closed real analytic subsets of $\Sigma$.  Since
$\Sigma$ is real analytic and connected, $\Sigma_j$ must either
coincide with $\Sigma$ or else be of zero measure in $\Sigma$ for
$j=1,2$. However, not both $\Sigma_1$ and $\Sigma_2$ can be of zero
measure in $\Sigma$ since $\Sigma_1\cup\Sigma_2=\Sigma$.  Hence,
$\Sigma_1=\Sigma$ or $\Sigma_2=\Sigma$, which completes the proof.
\end{proof}

It is worth making some remarks about Theorem \ref{ruledthm}. Note
that (i) and (ii) are \emph{linear} conditions on $\psi$ and, by the
remarks made above, \eq{assruled4} is the condition which makes the
asymptotic cone $M_0$ associative. So, if we start with an
associative two-sided cone $M_0$ defined by a map $\phi$, then
$\phi$ and a function $\psi$ satisfying (i) or (ii) will define an
r-oriented ruled associative 3-fold $M$ with asymptotic cone $M_0$.
We also note that conditions (i) and (ii) are unchanged if $\phi$ is
fixed and satisfies \eq{assruled4}, but $\psi$ is replaced by
$\psi+\tilde{f}\phi$ where $\tilde{f}$ is a real analytic function.
We can thus always locally transform $\psi$ such that $f$ in
condition (i) is zero.

\subsection{Evolution equations for ruled associative 3-folds}

Our first result follows \cite[Proposition 5.2]{Joy5}. Here we make
the definition that a function is real analytic on a compact
interval $I$ in $\R$ if it extends to a real analytic function on an
open set containing $I$.

\begin{thm}
\label{evolveruled} Let $I$ be a compact interval in $\R$, let $s$
be a coordinate on $I$, and let $\phi_0:I\rightarrow
\mathcal{S}^6$ and $\psi_0:I\rightarrow\R^7$ be real analytic
maps.
 Then there exist $\epsilon>0$ and unique real analytic maps $\phi:I\times
(-\epsilon,\epsilon)\rightarrow \mathcal{S}^6$ and $\psi:I\times
(-\epsilon,\epsilon) \rightarrow\R^7$ satisfying
$\phi(s,0)=\phi_0(s)$, $\psi(s,0)=\psi_0 (s)$ for all $s\in I$ and
\begin{equation}
\label{evolruled}
 \frac{\partial \phi}{\partial t} = \phi\times\frac{\partial \phi}{\partial s}\,,
\hspace{20pt} \frac{\partial \psi}{\partial t} =
\phi\times\frac{\partial \psi}{\partial s}\,,
\end{equation}
\noindent where $t$ is a coordinate on $(-\epsilon,\epsilon)$.  Let $M$ be
defined by:
\begin{equation*}
\begin{split}
M = \{r\phi(s,t)+\psi(s,t):r\in\R,\hspace{1pt} s\in I,\hspace{1pt} t\in(
-\epsilon,\epsilon)\}.
\end{split}
\end{equation*}
\noindent Then $M$ is an r-oriented ruled associative 3-fold in $\R^7$.
\end{thm}

\begin{proof}
Since $I$ is compact and $\phi_0,\psi_0$ are real analytic, we may
use the \emph{Cauchy--Kowalevsky Theorem} \cite[p.234]{Racke} to
give us functions $\phi:I\times(-\epsilon,\epsilon)\rightarrow\R^7$
and $\psi:I\times (-\epsilon,\epsilon)\rightarrow\R^7$ satisfying
the initial conditions and \eq{evolruled}.  It is clear that
$\frac{\partial}{\partial t}g(\phi,
\phi)=2g(\phi,\frac{\partial\phi}{\partial t})=0$, since
$\frac{\partial\phi}{\partial t}$ is defined by a cross product
involving $\phi$ and hence is orthogonal to $\phi$.  We may deduce
that $|\phi|$ is independent of $t$ and is therefore one, so that
$\phi$ maps to $\mathcal{S}^6$. We conclude that $M$ is an
r-oriented ruled associative 3-fold using (i) of Theorem
\ref{ruledthm}.
\end{proof}

\noindent Theorem \ref{evolveruled} shows that \eq{evolruled} can
be considered as evolution equations for maps $\phi,\psi$
satisfying (i) of Theorem \ref{ruledthm}.  We now show that
condition (ii) of Theorem \ref{ruledthm} does not produce any
interesting ruled associative 3-folds.  We make the definition
that two rulings $(\Sigma,\pi)$ and $(\tilde{\Sigma},\tilde{\pi})$
are \emph{distinct} if the families of affine straight lines
$\mathcal{F}_{\Sigma}=\{\pi^{-1}(\sigma):\sigma\in\Sigma\}$ and
$\mathcal{F}_{\tilde{\Sigma}}=\{\tilde{\pi}^{-1}(\tilde{\sigma}):
\tilde{\sigma}\in\tilde{\Sigma}\}$ are different.

\begin{prop}
\label{planar} Any r-oriented ruled associative 3-fold
$(M,\Sigma,\pi)$ satisfying condition (ii) but not (i) of Theorem
\ref{ruledthm} is locally isomorphic to an affine associative
3-plane in $\R^7$.
\end{prop}

\begin{proof}
By Theorem \ref{thm1}, $M$ is real analytic wherever it is
nonsingular and so we can take $(\Sigma,\pi)$ to be locally real
analytic.  Let $I=[0,1]$, let $\gamma:I\rightarrow\Sigma$ be a
real analytic curve in $\Sigma$ and let $\phi,\psi$ be the
functions defining $M$.
 Then we can use Theorem \ref{evolveruled} with initial conditions
$\phi_0=\phi(\gamma(s))$ and $\psi_0=\psi(\gamma(s))$ to give us
functions $\tilde{\phi},\tilde{\psi}$, which define an r-oriented
ruled associative 3-fold $\tilde{M}$ satisfying (i) of Theorem
\ref{ruledthm}. However, $M$ and $\tilde{M}$ coincide in the real
analytic\ 2-fold $\pi^{-1}(\gamma(I))$, and hence, by Theorem
\ref{thm2}, they must be locally equal.  We conclude that $M$
locally admits a ruling $(\tilde{\Sigma},\tilde{\pi})$ satisfying
(i) of Theorem \ref{ruledthm}, which must therefore be distinct
from the ruling $(\Sigma,\pi)$.


The families of affine straight lines $\mathcal{F}_{\Sigma}$ and
$\mathcal{F}_{\tilde{\Sigma}}$, using the notation above, coincide
in the family of affine straight lines defined by points on
$\gamma$, denoted $\mathcal{F}_{\gamma}$.  Using local real
analyticity of the families, either $\mathcal{F}_{\Sigma}$ is
equal to $\mathcal{F}_{\tilde{\Sigma}}$ locally or they only meet
in $\mathcal{F}_{\gamma}$ locally. The former possibility is
excluded because the rulings $(\Sigma,\pi)$ and
$(\tilde{\Sigma},\tilde{\pi})$ are distinct and thus the latter is
true.

Let $\gamma_1$ and $\gamma_2$ be distinct real analytic curves
near $\gamma$ in $\Sigma$ defining rulings $(\Sigma_1,\pi_1)$ and
$(\Sigma_2,\pi_2)$, respectively, as above.  Then
$\mathcal{F}_{\Sigma}\cap\mathcal{F}_{\Sigma_j}$ is locally equal
to $\mathcal{F}_{\gamma_j}$ for $j=1,2$.  Hence,
$(\Sigma_1,\pi_1)$ and $(\Sigma_2,\pi_2)$ are not distinct (that
is, $\mathcal{F}_{\Sigma_1}=\mathcal{F}_{\Sigma_2}$) if and only
if $\mathcal{F}_{\gamma_1}=\mathcal{F}_{\gamma_2}$, which implies
that $\gamma_1=\gamma_2$.  Therefore, distinct curves near
$\gamma$ in $\Sigma$ produce different rulings of $M$ and thus $M$
has infinitely many rulings.

Suppose that $\{\gamma_t:t\in\R\}$ is a one parameter family of
distinct curves
 near $\gamma$ in $\Sigma$, with
$\gamma_0=\gamma$. Each curve in the family defines a distinct
ruling $(\Sigma_t,\pi_t)$, and hence there exists $p\in M$ with
$M$ nonsingular at $p$ such that $L_t=\pi_t^{-1}(\pi_t(p))$ is not
constant as a line in $\R^7$. We therefore get a one parameter
family of lines $L_t$ in $M$ through $p$ with $\frac{dL_t}{dt}\neq
0$ at some point, i.e. such that $L_t$ changes nontrivially. We
have thus constructed a real analytic one-dimensional family of
lines $\{L_t:t\in\R\}$ whose total space is a real analytic 2-fold
$N$ contained in $M$.  Moreover, every line in $M$ through $p$ is
a line in the affine associative 3-plane $p+T_pM$, and so $N$ is
contained in $p+T_pM$. Then, since $N$ has nonsingular points in
the intersection between $M$ and $p+T_pM$, Theorem \ref{thm2}
shows that $M$ and $p+T_pM$ coincide on a connected component of
$M$. Hence, $M$ is planar, i.e. $M$ is locally isomorphic to an
affine associative 3-plane in $\R^7$.
\end{proof}

We now state our main result on ruled associative 3-folds, which
follows from the results in this section.

\begin{thm}
\label{ruledmainresult} Let $(M,\Sigma,\pi)$ be a non-planar,
r-oriented, ruled associative\ 3-fold in $\R^7$.  Then there exist
real analytic maps $\phi:\Sigma\rightarrow \mathcal{S}^6$ and
$\psi:\Sigma\rightarrow\R^7$ such that:
\begin{align}
M & =  \{r\phi(\sigma)+\psi(\sigma):r\in\R,\hspace{2pt} \sigma\in\Sigma\}, \nonumber\\
\label{ruledmaineqphi}\frac{\partial\phi}{\partial t} & =  \phi\times\frac{\partial\phi}{\partial s}\,,\\
\label{ruledmaineq} \frac{\partial\psi}{\partial t} &=
\phi\times\frac{\partial\psi}{\partial s}+f\phi,
\end{align}
\noindent where $(s,t)$ are oriented conformal coordinates on $\Sigma$ and
$f:\Sigma\rightarrow\R$ is some real analytic function.

Conversely, suppose $\phi:\Sigma\rightarrow\ \mathcal{S}^6$ and
$\psi:\Sigma\rightarrow\R^7$ are real analytic maps satisfying
\eq{ruledmaineqphi} and \eq{ruledmaineq} on a connected real
analytic 2-fold $\Sigma$. If $M$ is defined as above, then $M$ is
an r-oriented ruled associative 3-fold wherever it is nonsingular.
\end{thm}

\subsection{Holomorphic vector fields}

We now follow \cite[$\S 6$]{Joy5} and use a \emph{holomorphic
vector field} on a Riemann surface $\Sigma$ to construct ruled
associative 3-folds.

\begin{prop}
\label{Lieass} Let $M_0$ be an r-oriented two-sided associative
cone in $\R^7$.  We can then write $M_0$ in the form \eq{ruled2}
for a real analytic map $\phi:\Sigma \rightarrow \mathcal{S}^6$,
where $\Sigma$ is a Riemann surface.  Let $w$ be a holomorphic
vector field on $\Sigma$ and define a map $\psi:\Sigma\rightarrow
\R^7$ by $\psi=\mathcal{L}_w\phi$, where $\mathcal{L}_w$ is the
Lie derivative with respect to $w$.  If we define $M$ by equation
\eq{ruled1} then $M$ is an r-oriented ruled associative 3-fold in
$\R^7$ with asymptotic cone $M_0$.
\end{prop}

\begin{proof}
We need only consider the case where $w$ is not identically zero
since this is trivial.  Then $w$ has only isolated zeros and,
since the fact that $M$ is associative is a closed condition on
the nonsingular part of $M$, it is sufficient to prove that
\eq{ruledmaineq} holds at any point $\sigma\in\Sigma$ such that
$w(\sigma)\neq 0$.  Suppose $\sigma$ is such a point. Then, since
$w$ is a holomorphic vector field, there exists an open set in
$\Sigma$ containing $\sigma$ on which oriented conformal
coordinates $(s,t)$ may be chosen such that
$w=\frac{\partial}{\partial s}$. Hence,
$\psi=\frac{\partial\phi}{\partial s}$ in a neighbourhood of
$\sigma$, and differentiating \eq{ruledmaineqphi} gives:
\begin{equation*}
\begin{split}
\frac{\partial^2\phi}{\partial s\partial t}=\frac{\partial\phi}{\partial s}
\times\frac{\partial\phi}{\partial s} + \phi\times\frac{\partial^2\phi}{\partial s^2}.\end{split}
\end{equation*}
\noindent Interchanging the order of the partial derivatives on
the left-hand side and noting that the cross product is
alternating, we have that
\begin{equation*}
\begin{split}
\frac{\partial\psi}{\partial t}=\frac{\partial^2\phi}{\partial s \partial t}=
\phi\times\frac{\partial\psi}{\partial s}.
\end{split}
\end{equation*}
\noindent The result follows from Theorem \ref{ruledmainresult}.
\end{proof}

Having proved a result which enables us to construct ruled
associative 3-folds given an associative cone on a Riemann surface
$\Sigma$, we consider which choices for $\Sigma$ will produce
interesting examples.  The only nontrivial vector spaces for
holomorphic vector fields on a compact connected Riemann surface
occur for genus zero or one.  We therefore focus our attention
upon the cases where we take $\Sigma$ to be $\mathcal{S}^2$ or
$T^2$. The space of holomorphic vector fields on $\mathcal{S}^2$
is 6-dimensional, and on $T^2$ it is 2-dimensional.  In the SL
case, any SL cone on $\mathcal{S}^2$ has to be an SL 3-plane
\cite[Theorem B]{Haskins}; Bryant \cite[$\S$4]{Bryant} shows that
this is not true in the associative case and that, in fact, there
are many nontrivial associative cones on $\mathcal{S}^2$.

\begin{thm}\label{holoasym}
Let $M_0$ be an r-oriented, two-sided, associative cone on a
Riemann surface $\Sigma \cong \mathcal{S}^2$ (or $T^2$) with
associated real analytic map $\phi:\Sigma\rightarrow
\mathcal{S}^6$ as in \eq{ruled2}.  Then there exists a
6-dimensional (or 2-dimensional) family of distinct r-oriented
ruled associative 3-folds with asymptotic cone $M_0$, which are
asymptotically conical to $M_0$ with order $O(r^{-1})$.
\end{thm}

\begin{proof}
If $(s,t)$ are oriented conformal coordinates on $\Sigma$, we may write holomorphic vector fields
on $\Sigma$ in the form:
\begin{equation}
\begin{split}
\label{holo} w = u(s,t)\frac{\partial}{\partial s} +
v(s,t)\frac{\partial}{\partial t}\,,
\end{split}
\end{equation}
\noindent where $u,v:\R^2\rightarrow \R$ satisfy the Cauchy--Riemann
equations.  For each holomorphic vector field $w$, as written in
\eq{holo}, define a 3-fold $M_w$ by:
\begin{equation*}
\begin{split}
M_w=\left\{r\phi(s,t)+u(s,t)\frac{\partial\phi}{\partial s}(s,t)+v(s,t)\frac{\partial\phi}
{\partial t}(s,t):r\in\R,(s,t)\in\Sigma\right\}.
\end{split}
\end{equation*}
By Proposition \ref{Lieass}, $M_w$ is an r-oriented ruled
associative 3-fold with asymptotic cone $M_0$, and it is clear
that each holomorphic vector field $w$ will give a distinct
3-fold.

We now construct a diffeomorphism $\Phi$ as in Definition
\ref{asym} satisfying \eq{asymconds} for $\alpha=-1$.  Let $R>0$,
$w$ be a holomorphic vector field as in \eq{holo}, and let
$\bar{B}_R$ denote the closed ball of radius $R$ in $\R^7$.
Define $\Phi:M_0\setminus\bar{B}_R\rightarrow M_w$ by:
\begin{equation*}
\begin{split}
\Phi(r\phi(s,t))=r\phi\left(s-\frac{u}{r},t-\frac{v}{r}\right)\!+\!
u\frac{\partial\phi}{\partial
s}\left(s-\frac{u}{r},t-\frac{v}{r}\right)\!+\!
v\frac{\partial\phi}{\partial
t}\left(s-\frac{u}{r},t-\frac{v}{r}\right),
\end{split}
\end{equation*}
\noindent where $|r|>R$.  Clearly, $\Phi$ is a well-defined map
with image in $M_w\setminus K$ for some compact subset $K$ of
$M_w$. Note that, by choosing $R$ sufficiently large, we can
expand the various terms defining $\Phi$ in powers of $r^{-1}$ as
follows:
\begin{align*}
\phi\left(\!s\!-\frac{u(s,t)}{r},t\!-\frac{v(s,t)}{r}\right)\! & =
\phi(s,t)\!-\frac{u(s,t)}{r} \frac{\partial\phi}{\partial
s}(s,t)\!-\frac{v(s,t)}{r}\frac{\partial\phi}{\partial t}
+ O(r^{-2}),\\
\frac{\partial\phi}{\partial
s}\left(\!s\!-\frac{u(s,t)}{r},t\!-\frac{v(s,t)}{r}\right)\! &=
\frac{\partial\phi}{\partial s}(s,t) + O(r^{-1}), \\
\frac{\partial\phi}{\partial
t}\left(\!s\!-\frac{u(s,t)}{r},t\!-\frac{v(s,t)}{r}\right)\! &=
\frac{\partial\phi}{\partial t}(s,t) + O(r^{-1}).
\end{align*}
\noindent We deduce that
\begin{equation*}
\begin{split}
|\Phi(r\phi(s,t))-r\phi(s,t)| = O(r^{-1}),
\end{split}
\end{equation*}
\noindent and the other conditions in \eq{asymconds} can be
derived similarly. We conclude that $M_w$ is asymptotically
conical to $M_0$ with order $O(r^{-1})$.
\end{proof}

\noindent There are many examples of associative cones over $T^2$
given by the SL tori constructed by Haskins \cite{Haskins}, Joyce
\cite{Joy2} and McIntosh \cite{McIntosh} and others. However, by
Theorem \ref{SLconverge}, applying Theorem \ref{holoasym} to them
will only produce ruled SL 3-folds and the result reduces to
\cite[Theorem 6.3]{Joy5}.


\end{document}